\documentclass[article,11pt]{elsarticle}

\usepackage{amsmath, amsthm, amssymb}
\usepackage{graphicx}
\usepackage{psfrag}
\usepackage{epsfig}
\usepackage{epstopdf}

\newtheorem{theorem}{Theorem}

\newcommand{{\Reals}}{\mathbb{R}}

\begin{document}

\begin{frontmatter}

\title{A partitioned scheme for fluid-composite structure interaction problems \tnoteref{label1}}

 \author[lab1]{M. Buka\v{c}, S. \v{C}ani\'{c}, B. Muha}




\begin{abstract}
We present a loosely-coupled partitioned scheme for a benchmark problem in fluid-composite structure interaction.
The benchmark problem proposed here consists of an incompressible, viscous fluid interacting with a composite structure
that consists of two layers: a thin elastic layer which is in contact with the fluid, and 
a thick elastic layer which sits on top of the thin layer. The motivation comes from fluid-structure interaction (FSI) in hemodyamics. 
The equations of linear elasticity are used to model the thick structural layer, while the Koiter member/shell equations
are used to model the thin structural layer which serves as fluid-structure interface with mass.
An efficient, modular, operator-splitting scheme is proposed to simulate solutions to the coupled, nonlinear FSI problem. 
The operator splitting scheme separates the elastodynamics structure problem,
from a fluid problem in which the thin structure inertia is included as a Robin-type boundary condition to achieve unconditional stability,
without requiring any sub-iterations within time-steps. 
An energy estimate associated with unconditional stability is derived
for the fully nonlinear FSI problem defined on  moving domains. 
Two instructive numerical examples are presented to test the performance of the scheme, where it is shown numerically, that the scheme is
at least first-order accurate in time. The second example reveals a new phenomenon in FSI problems: the presence of a thin fluid-structure interface
with mass regularizes solutions to the full FSI problem. 
\end{abstract}

\begin{keyword}
fluid-structure interaction \sep
composite structure \sep
partitioned scheme \sep
blood flow

\end{keyword}

\end{frontmatter}

\section{Introduction}

\subsection{Motivation}

Composite materials are materials made from two or more constituent materials with different physical properties 
that when combined, produce a material with characteristics different from the individual components \cite{wikipedia}.
The new material may be preferred for many reasons: common examples include materials which are stronger, 
lighter or less expensive when compared to traditional materials.
Composite materials have been used in virtually all areas of engineering. Examples
include marine structures such as boats, ships, and offshore structures,
in aerospace industry such as aircraft design, in the design of wind and marine turbines, and in the design of sports' 
equipment such as skis and tennis raquets.
They also appear everywhere in nature. Examples include bones, horns, wood, and biological tissues such as  blood vessels of major human arteries.
These materials are exposed to a wide spectrum of dynamic loads. 
In marine structures, aerospace industry, and in many biological constructs, weather natural or man-made, 
the dynamic loading comes from the surrounding fluid such as water, air, or blood flow.
Understanding the interaction between composite materials and the surrounding fluid is
important for the understanding of the normal function of the coupled fluid-structure system, as well as
the detection of damage and/or pathologies in their function. Prevention of catastrophic events in engineered constructs,
or the design of medical treatments to prevent further biological tissue damage,
is aided by computational and experimental studies of fluid-composite structure interaction.

In problems in which the fluid and structure have comparable densities, which is the case with most biological tissues (e.g.,
arterial walls and blood), or in problems in which the structure density is less than the density of water, which is the case in several sandwich 
composite materials immersed in water \cite{KwonOwens}, the exchange in energy between the fluid and structure is so significant that classical Dirichlet-Neumann
loosely-coupled partitioned schemes for numerical simulation of these kinds of problems are unconditionally unstable 
\cite{causin2005added,Michler}. This is due to the so called added mass effect, also known as virtual mass, or hydrodynamic mass.
As the structure vibrates, its mass is increased by the mass of the surrounding fluid, consequently decreasing its natural frequency. 
This is due to the increase of total kinetic energy of the structure and fluid from the addition of the kinetic energy of the fluid,
which significantly impacts the structure vibrations when the structure is ``light'' relative to the fluid. 
A comprehensive study of these problems remains to be
a challenge due to their significant energy exchange leading to strong nonlinearities in the 
coupled FSI problem, and due to the intrinsic multi-physics nature of the problem.

\begin{figure}[ht]
\centering{
\includegraphics[scale=0.25]{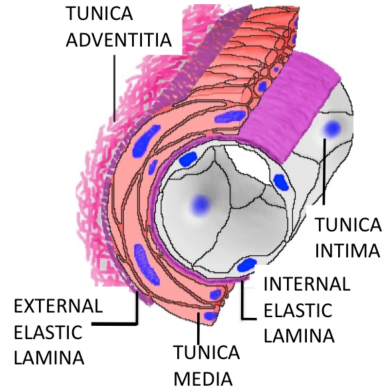}
}
\caption{Arterial wall layers.}
\label{fig:aretria}
\end{figure}

The motivation for this paper comes from fluid-structure interaction between blood flow and cardiovascular tissue. 
Arterial walls of major human arteries are composite materials consisting of three main layers: the intima, media and adventitia,
separated by thin elastic laminae, see Figure~\ref{fig:aretria}.
The intima is the innermost, thin layer, which is mainly composed of endothelial cells. The media is the middle 
layer and it is mainly composed of elongated smooth muscle cells, and also elastin and collagen. 
The adventitia is the outermost layer, and it  is manly composed of collagen fibrils, elastic sheets and elastic fibrils.
The media gives rise to the majority of the vessel's elastic (viscoelastic) behavior.

Arterial walls interact with blood flow. This interaction plays a crucial role in the normal physiology and pathophysiology of the human
cardiovascular system. During systole, when the left ventricle of the heart
contracts and squeezes blood to the aorta, the kinetic energy imparted by blood to the aortic walls
stretches the walls outward. 
Once the artery is 
distended with blood, energy stored by stretching elastic fibers is released through elastic recoil.  Elastic recoil takes place during the diastolic 
part of cardiac cycle, when the heart's left ventricle relaxes and gets refilled with blood.  
During that time the elastic recoil of arteries helps propel blood to the far most
parts of the cardiovascular system. 

In medium-to-large human arteries, such as the aorta or coronary arteries, blood can be modeled as an incompressible,
viscous, Newtonian fluid. Due to the fact that density of blood is comparable to that of arterial wall tissue, 
the corresponding fluid-structure interaction problem suffers from the added mass effect, discussed above.
Except for the Immersed Boundary Method approaches,
to this date, there have been no FSI numerical solvers based on the Arbitrary Lagrangian-Eulerian approaches,
that would take into account the multi-layered structure of arterial walls,
and the full nonlinear (two-way) coupling with blood flow.
In this manuscript we take a first step in this direction by proposing to study 
a benchmark problem in fluid-multi-layered-structure interaction
in which the structure consists of two layers, a thin and a thick layer, 
where the thin layer serves as fluid-structure interface with mass.
The thin layer is modeled by the cylindrical Koiter membrane equations,
while the thick layer by the equations of linear elasticity. 
The proposed problem is a nonlinear moving-boundary problem of parabolic-hyperbolic 
type. 

Our computational method is based on an operator-splitting approach, introduced
in \cite{BorSunMulti} to study the existence of solutions to an associated fluid-multi-layered structure interaction problem. 
The particular Lie operator splitting used in \cite{BorSunMulti} gives rise to a modular scheme that deals efficiently and elegantly with the
added mass effect problem. In this scheme, the time discretization is performed via Lie operator splitting.
At each time step, the full FSI problem is split into a fluid and a structure sub-problem. 
To deal with the motion of the fluid domain, the Arbitrary Lagrangian Eulerian (ALE) approach is adopted.
To achieve stability and convergence, the splitting
is performed in a special way in which the fluid sub-problem 
includes the thin structure inertia  via a  "Robin-type" boundary condition.  
The fact that the thin structure inertia is included implicitly in the fluid sub-problem, gives rise to the
appropriate energy estimates for the approximate solutions, independent of the size of the time discretization,
that provide unconditional stability of the underlying scheme \cite{BorSunMulti}. 

In this manuscript, the proposed numerical scheme is implemented to solve two instructive numerical examples.
The first example shows that the scheme recovers an exact solution
to a simplified FSI problem. 
The second example considers a full FSI problem with a thin and thick structural layer. 
Since there are no numerical results in literature on FSI problems with multiple structural layers against which the numerical solution
could be tested, the following strategy was used to test the solver.
A sequence of solutions to  FSI problems with a thin and thick structure was considered 
in which the thickness of the thin structure $h$ converges to zero, while the combined thickness of the composite structure remains constant.
Using analytical methods, it is shown in this manuscript that the limiting solution as $h\to 0$
solves the FSI problem with a single, thick structural layer. 
The solution of the limiting problem is then numerically tested against the solution of the FSI problem with the  thick structural layer,
which is obtained using a different solver, published in \cite{thick}, and previously tested against a monolithic solver. 
It is then shown that the two solutions,
obtained with two different solvers, are in good agreement.
In this example it is also shown numerically that the proposed operator-splitting scheme is first-order accurate in time.

Although the benchmark problem and numerical examples in this manuscript are given in 2D, 
there is nothing in the methodology that depends on the spatial dimension. Therefore, the methodology presented
in this manuscript can be applied to 3D FSI problems. An application of the splitting to a problem in  3D can be found in \cite{BorSun3D},
where a FSI with a single thin structure was considered. 

We also emphasize that the same (equivalent) splitting strategy can be applied to problems in which the composite structure is viscoelastic. 
Such a splitting was presented, for example, in \cite{Martina_paper1}, where a single, thin viscoelastic Koiter shell was coupled
to the motion of the fluid. See also \cite{thick} for a corresponding FSI problem involving a thick viscoelastic structure.  
Since structural viscosity further stabilizes the underlying FSI problem, in this manuscript we consider the most ``difficult'' case, 
namely the case when the composite structure is purely elastic.

The work presented in this manuscript reveals a new phenomenon which is related to the smoothing effects provided by the 
purely elastic fluid-structure interface {\sl with mass}.
Our simulations show that the presence of a thin fluid-structure interface with mass regularizes solutions of the entire FSI problem. 
Our results in the companion paper on the analysis (existence) of solutions to fluid-multi-layered structure interaction problem \cite{BorSunMulti} 
show that the {\sl inertia} of the thin fluid-structure interface with mass provides higher regularity of the coupled problem. 
This is reminiscent of the results by Hansen and Zuazua \cite{HansenZuazua} in which the presence of point mass
at the interface between two linearly elastic strings with solutions in asymmetric spaces (different
regularity on each side) allowed the proof of well-posedness due to the regularizing effects by the
point mass. 
Based on the computational results presented in example 2 of this manuscript,
we conclude that the fluid-structure interface inertia regularizes the motion of the fluid-structure interface,
and acts as a regularizing mechanism for  the entire solution of this FSI problem.

\subsection{Description of the fluid-structure interaction problem}\label{description}

\begin{figure}[ht]
\centering{
\includegraphics[scale=0.6]{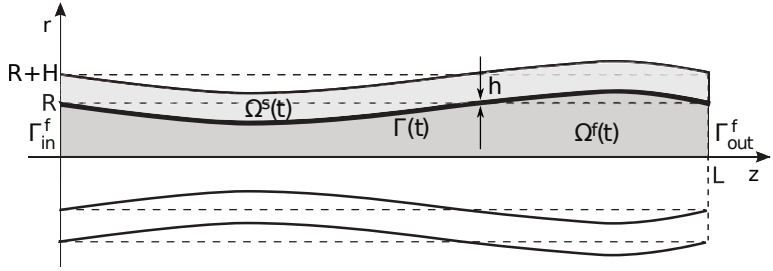}
}
\caption{Deformed domains $\Omega^f(t) \cup \Omega^s(t).$}
\label{fig:domain}
\end{figure}

We consider the flow of an incompressible, viscous fluid in a two-dimensional channel of reference length $L$, and reference width $2R$,
see Figure~\ref{fig:domain}.
The channel is bounded by a deformable wall, which is composed of two layers: a thin elastic layer with thickness $h$, and a thick elastic layer with thickness $H$. The thin layer serves as fluid-structure interface with mass. 

We are interested in simulating a pressure-driven flow through the deformable 2D channel with two-way coupling 
between the fluid and structure. The structure deforms because of the fluid loading onto the structure, 
while at the same time, the fluid flow and the fluid domain are affected by the 
elastodynamics of the structure. Deformation of the structure will be given with respect of a reference configuration, i.e., 
in Lagrangian framework,
where the reference structural domain is assumed to be a uniform cylindrical annulus with constant radius and thickness. 
Without loss of generality, we consider only the upper half of the fluid domain supplemented by a symmetry condition at the axis of symmetry.
Thus, the reference fluid and structure domains in our problem are given, respectively, by
\begin{eqnarray*}
\hat{\Omega}^f &:=& \{(z,r) | 0<z<L, 0<r<R \}, \\
\hat{\Omega}^s &:=& \{(z,r) | 0<z<L, R<r<R+H \}.
\end{eqnarray*}
Here $z$ and $r$ denote the horizontal and vertical Cartesian coordinates, respectively (see Figure~\ref{fig:domain}).
The reference configuration of the fluid-structure interface is given by
\begin{eqnarray*}
\hat{\Gamma} &:=& \{(z,R) | 0<z<L \}.
\end{eqnarray*}

The flow of an incompressible, viscous fluid is modeled by  the Navier-Stokes equations:
\begin{align}\label{NS1}
& \rho_f \bigg( \frac{\partial \boldsymbol{v}}{\partial t}+ \boldsymbol{v} \cdot \nabla \boldsymbol{v} \bigg)  = \nabla \cdot \boldsymbol\sigma^f(\boldsymbol v, p) &  \textrm{in}\; \Omega^f(t)\times(0,T), \\
 \label{NS2}
&\nabla \cdot \boldsymbol{v} = 0 & \textrm{in}\; \Omega^f(t)\times(0,T),
\end{align}
where $\boldsymbol{v}=(v_z,v_r)$ is the fluid velocity, $p$ is the fluid pressure, $\rho_f$ is fluid density, and $\boldsymbol\sigma^f$ is the fluid stress tensor. 
For a Newtonian fluid the stress tensor is given by $\boldsymbol\sigma^f(\boldsymbol v,p) = -p \boldsymbol{I} + 2 \mu_f \boldsymbol{D}(\boldsymbol{v}),$ where
$\mu_f$ is the fluid viscosity and  $\boldsymbol{D}(\boldsymbol{v}) = (\nabla \boldsymbol{v}+(\nabla \boldsymbol{v})^{\tau})/2$ is the rate-of-strain tensor.

The inlet to the fluid domain will be denoted by $\Gamma^f_{in}= \{0\} \times (0,R)$, 
and the outlet by $\Gamma^f_{out} = \{L\}\times (0,R)$.
The flow will be driven by the inlet and outlet ``pressure'' data.
Two sets of inlet and outlet boundary conditions will be considered. 
\begin{enumerate}
\item[ ] {\bf BC1.} The {\sl normal stress} inlet and outlet data:
\begin{equation}
\boldsymbol {\sigma n}^f_{in/out} = -p_{in/out}(t) \boldsymbol{n}^f_{in/out} \  \textrm{on} \; \Gamma^f_{in/out} \times (0,T), \label{inlet} 
\end{equation}
where $\boldsymbol{n}^f_{in}$ and $\boldsymbol{n}^f_{out}$ are the corresponding outward unit normals.
Even though not physiologically optimal, these boundary conditions are common in blood flow simulations~\cite{badia2008fluid, miller2005computational, nobile2001numerical}.
\item[ ] {\bf BC2.}  The {\sl dynamic pressure} inlet and outlet data:
\begin{equation}\label{PD}
\left. \begin{array}{rcl}
\displaystyle{p+\frac{\rho_f}{2}|u|^2}&=&p_{in/out}(t) \\
u_r &=& 0\\
\end{array}
\right\} \quad {\rm on}\ \Gamma^f_{in/out} \times (0,T),
\end{equation}
where $p_{in/out}\in L^{2}_{loc}(0,\infty)$ are given. 
Therefore, the fluid flow is driven by a prescribed dynamic pressure drop, 
and the flow enters and leaves the fluid domain orthogonally to the inlet and outlet boundary. 
This set of boundary conditions will be used in energy estimates and in the proof of stability 
of the proposed scheme. 
\end{enumerate}

At the bottom fluid boundary denoted by $\Gamma_b$, defined by $(0,L)\times\{ 0 \}$,
 the following symmetry conditions are prescribed:
\begin{equation}\label{symmetry_condition}
 \frac{\partial v_z}{\partial r}(z,0,t) = 0, \quad v_r(z,0,t) = 0 \quad \textrm{on} \; (0,L) \times (0,T).
\end{equation}

The lateral fluid boundary is bounded by a deformable, composite structure. 
The thin layer, which is in contact with the fluid, 
is modeled by the linearly elastic Koiter membrane model:
\begin{align}
& \rho_{m} h \frac{\partial^2 \hat\eta_z}{\partial t^2}-C_2 \frac{\partial \hat\eta_r}{\partial \hat{z}}-C_1 \frac{\partial^2 \hat\eta_z}{\partial\hat {z}^2}= \hat {f}_z  & \textrm{on} \; \hat{\Gamma} \times (0,T), \label{structure1}  \\
& \rho_{m} h \frac{\partial^2 \hat\eta_r}{\partial t^2}+C_0 \hat\eta_r +C_2 \frac{\partial \hat\eta_z}{\partial \hat{z}}   = \hat {f}_r & \textrm{on} \; \hat{\Gamma}\times (0,T), \label{structure2}
\end{align}
where $\hat{\boldsymbol \eta} (\hat{z},t)= (\hat\eta_z(\hat{z},t), \hat\eta_r(\hat{z},t))$ denotes the axial and radial components of displacement, 
$\hat{\boldsymbol f} = (\hat {f}_z, \hat {f}_r)$ is the surface density of the force applied to the shell, 
$\rho_{m}$ denotes membrane density, and 
{\small{
\begin{equation}\label{coefficientsC}
 C_0 = \frac{h}{R^2} \big(\frac{2 \mu_{m} \lambda_{m}}{\lambda_{m}+2\mu_{m}}+2\mu_{m}\big), \;\;  C_1 = h \big(\frac{2 \mu_{m} \lambda_{m}}{\lambda_{m}+2\mu_{m}}+2\mu_{m}\big), \;\; C_2 =\frac{h}{R} \frac{2 \mu_{m} \lambda_{m}}{\lambda_{m}+2\mu_{m}},
\end{equation}}}
where  $\mu_{m}$ and $\lambda_{m}$ are the Lam\'e constants. In terms of the Young's modulus of elasticity and Poisson ratio,
coefficients $C_0$, $C_1$, and $C_2$ can be written as
 \begin{equation*}
 C_0 = \frac{hE}{R^2(1-\sigma^2)}, \;\;  C_1 = \frac{hE}{1-\sigma^2}, \;\; C_2 =\frac{hE\sigma}{R(1-\sigma^2)}.
\end{equation*}

The thick layer of the wall will be modeled by the equations of linear elasticity with an added extra term
$\gamma \hat{\bf U}$, where $\hat{\bf U}$ denotes the thick structure displacement,
 to account for circumferential strain whose effects are lost in the transition from  3D to  2D \cite{badia2008splitting,badia,Ma,barker2010scalable}.
This term corresponds to the non-differentiated term in the Koiter membrane equations \eqref{structure2} containing the coefficient $C_0$,
which appears in these equations due to the cylindrical geometry of the domain.
If the structure is not fixed at the end points, this term helps keep the top and bottom portions of the structure domain together.
 The model reads:
\begin{equation}
\rho_{s} \frac{\partial^2 \hat{\boldsymbol U}}{\partial t^2} + \gamma \hat{\boldsymbol U}= \nabla \cdot \boldsymbol \sigma^s( \hat{\boldsymbol U}) \quad \textrm{in} \; \hat{\Omega}^s \times (0,T), \label{linela}
\end{equation}
where $\hat{\boldsymbol U} = (\hat{U}_z, \hat{U}_r)$ is the structure displacement and  $\rho_{s}$ is the structure density.  
We will be assuming that the thick structure behaves as a linearly elastic Saint-Venant Kirchhoff material, in which case 
the stress-displacement relationship is given by
$$
\boldsymbol \sigma^s( \hat{\boldsymbol U}) = 2 \mu_{s} \boldsymbol D(\hat{\boldsymbol U}) + \lambda_{s} (\nabla \cdot \hat{\boldsymbol U}) \boldsymbol I,
$$
where $\mu_{s}$ and $\lambda_{s}$ are the Lam\'e coefficients for the thick layer. 
In the simulations presented in Sections~\ref{sec:numerics1} and \ref{sec:numerics2}
the Lam\'{e} constants of the thin and thick structural layer will be assumed
to be the same:
$$
\lambda_m = \lambda_s, \ \mu_m = \mu_s.
$$

The thick structure will be assumed to be fixed at the inlet and outlet boundaries:
\begin{equation}\label{homostructure1}
 \hat{\boldsymbol U}(0, r, t) = \hat{\boldsymbol U}(L, r, t) = 0 \quad \textrm{on} \; [R, R+H]\times(0,T).
\end{equation}
Furthermore, at the external structure boundary $\hat{\Gamma}^s_{ext} = \{R+H\} \times (0,L)$ 
the structure will be exposed to zero external ambient pressure, while the axial displacement remains fixed:
\begin{align}\label{homostructure2}
& \boldsymbol n^s_{ext} \cdot \boldsymbol \sigma^s \boldsymbol { n}^s_{ext} =  0 & \textrm{on} \; \hat{\Gamma}^s_{ext} \times (0,T), \\
& \hat{U}_z =  0 & \textrm{on} \; \hat{\Gamma}^s_{ext} \times (0,T),
\end{align}
where $\boldsymbol n^s_{ext}$ is the outward unit normal vector on $\hat{\Gamma}^s_{ext}$.

Initially, the fluid and the structure are assumed to be at rest, with zero displacement from the reference configuration
\begin{equation}\label{initial}
\boldsymbol{v}=0, \quad \hat{\boldsymbol \eta}=0,
 \quad \frac{\partial \hat{\boldsymbol \eta}}{\partial t}=0, \quad \hat{\boldsymbol U} = 0, 
 \quad \frac{\partial \hat{\boldsymbol U}}{\partial t}=0,\quad {\rm at}\; t=0.
\end{equation}

The fluid and the composite structure are coupled via the following kinematic and dynamic boundary conditions~\cite{canic2005self,mikelic2007fluid}:
\begin{itemize}
\item \textbf{Kinematic coupling conditions} describe continuity of  velocity at the fluid-structure interface (no-slip)
\begin{equation}\label{kinematic}
 \boldsymbol{v}(\hat{z}+\hat{\eta}_z(\hat{z},t),R+\hat{\eta}_r(\hat{z},t),t)=\frac{\partial \hat{\boldsymbol \eta}}{\partial t}(\hat{z}, t) \quad \; \textrm{on} \; (0,L)\times (0,T),
\end{equation}  
and continuity of displacement
\begin{equation}\label{kinematic1}
\hat{\boldsymbol \eta}(\hat{z}, t)= \hat{\boldsymbol U}(\hat{z}, R, t) \quad \; \textrm{on} \; (0,L)\times (0,T);
\end{equation} 
\item \textbf{Dynamic coupling condition} describes the second Newton's law of motion of the fluid-structure interface, which is 
loaded by the normal stresses exerted by both the fluid and thick structure. The condition reads:
\begin{equation}
 \left(
 \begin{array}{l}
\displaystyle{ \rho_{m} h \frac{\partial^2 \hat\eta_z}{\partial t^2}-C_2 \frac{\partial \hat\eta_r}{\partial \hat{z}}-C_1 \frac{\partial^2 \hat\eta_z}{\partial\hat {z}^2}} \\
\displaystyle{\rho_{m} h \frac{\partial^2 \hat\eta_r}{\partial t^2}+C_0 \hat\eta_r +C_2 \frac{\partial \hat\eta_z}{\partial \hat{z}}  }
 \end{array}
 \right)
 =J\ \widehat{\boldsymbol{\sigma}^f\boldsymbol { n}^f|_{\Gamma(t)}} +  \boldsymbol \sigma^s \boldsymbol{n}^s|_{\hat{\Gamma}},
 \label{dynamic}
\end{equation}
on $(0,L)\times(0,T)$, where $J$ denotes the Jacobian of the transformation from Eulerian to Lagrangian coordinates, 
and $\widehat{\boldsymbol{\sigma}^f \boldsymbol {n}^f}|_{\Gamma(t)}$ denotes the normal fluid stress at the deformed fluid-structure interface,
evaluated on the reference configuration, namely at the points $(\hat{z}+\hat{\eta}_z(\hat{z},t),R+\hat{\eta}_r(\hat{z},t))$ for $\hat{z} \in (0,L)$
and $t \in (0,T)$.
Vector $\boldsymbol n^f$ is the outward unit normal to the deformed fluid domain, and $\boldsymbol n^s$ is the outward unit normal to the structure domain.
\end{itemize}

\subsection{The energy of the coupled problem}
The coupled fluid-composite structure interaction problem~\eqref{NS1}-\eqref{dynamic} with {\sl dynamic pressure} inlet and 
outlet data, satisfies the following energy equality:
{\small{
 \begin{eqnarray}
& \displaystyle\frac{d}{dt} \Bigg\{ \underbrace{ \frac{\rho_f}{2} ||\boldsymbol v||^2_{L^2(\Omega^f(t))}+ \frac{\rho_{m} h}{2} \bigg|\bigg|\frac{\partial \hat\eta_z}{\partial t}\bigg|\bigg|^2_{L^2(0,L)} 
+ \frac{\rho_{m} h}{2} \bigg|\bigg|\frac{\partial \hat\eta_r}{\partial t}\bigg|\bigg|^2_{L^2(0,L)} +  \frac{\rho_{s}}{2} \bigg|\bigg|\frac{\partial \hat{\boldsymbol U}}{\partial t} \bigg|\bigg|^2_{L^2(\hat{\Omega}^s)}}_{Kinetic \ Energy\ of\ the\ Coupled\ Problem}
 \nonumber\\
&
  \underbrace{+\displaystyle h \left[
  4 \mu_{m} \bigg|\bigg|\frac{\hat\eta_r}{R}\bigg|\bigg|^2_{L^2(0,L)} +
  4 \mu_{m} \bigg|\bigg|\frac{\partial \hat\eta_z}{\partial \hat{z}} \bigg|\bigg|^2_{L^2(0,L)}
  +\frac{4 \mu_{m} \lambda_{m}}{\lambda_{m}+2\mu_{m}} \bigg|\bigg|\frac{\partial \hat\eta_z}{\partial \hat{z}}+ \frac{\hat\eta_r}{R} \bigg|\bigg|^2_{L^2(0,L)}
 \right] }_{Membrane \ Elastic \ Energy\ } \nonumber \\
&  \underbrace{+ \frac{\gamma}{2}  ||\hat{\boldsymbol U}||^2_{L^2(\hat{\Omega}^s)}+\mu_{s}  ||\boldsymbol D(\hat{\boldsymbol U})||^2_{L^2(\hat{\Omega}^s)}+\frac{\lambda_{s}}{2} ||\nabla \cdot \hat{\boldsymbol  U}||^2_{L^2(\hat{\Omega}^s)} }_{Thick \ Structure \ Elastic \ Energy\ }
\Bigg\}
 \nonumber
\\
 & +  \underbrace{2 \mu_f ||\boldsymbol D(\boldsymbol v)||^2_{L^2(\Omega^f(t))}}_{Fluid \ Viscous \ Energy}
 = 
  \displaystyle{\int_0^R p_{in}(t) v_z|_{z=0} dr -   \int_0^R p_{out}(t) v_z|_{z=L} dr.}
 \label{FSIEnergy}
\end{eqnarray}
}}

The derivation of the energy equality above is similar to that in \cite{BorSunMulti} where the proof of existence of a weak solution
to the coupled problem was proved. 
The only difference with the energy obtained in \cite{BorSunMulti} is that the Koiter membrane model in the current manuscript
includes both the longitudinal and radial components of displacement. In this case, the total energy associated with the
Koiter membrane elastodynamics is given by
\begin{equation*}
\frac{d}{dt}\Big (\frac{\rho_{m} h}{2} \bigg|\bigg|\frac{\partial \hat\eta_z}{\partial t}\bigg|\bigg|^2_{L^2(0,L)}
+h \Big[
  4 \mu_{m} \bigg|\bigg|\frac{\hat\eta_r}{R}\bigg|\bigg|^2_{L^2(0,L)} +
  4 \mu_{m} \bigg|\bigg|\frac{\partial \hat\eta_z}{\partial \hat{z}} \bigg|\bigg|^2_{L^2(0,L)}
\end{equation*}  
\begin{equation}\label{KoiterEnergy}
+\frac{4 \mu_{m} \lambda_{m}}{\lambda_{m}+2\mu_{m}} \bigg|\bigg|\frac{\partial \hat\eta_z}{\partial \hat{z}}+ \frac{\hat\eta_r}{R} \bigg|\bigg|^2_{L^2(0,L)}
 \Big]  \Big )=\int_0^L \hat{\boldsymbol f} \cdot\frac{\partial \hat{\boldsymbol \eta}}{\partial t}  d\hat{z}.
\end{equation}
Equation \eqref{KoiterEnergy} is obtained after multiplying (using the dot-product) 
the Koiter membrane equations \eqref{structure1}, \eqref{structure2} 
by the structure velocity ${\partial\hat{\boldsymbol\eta}}/{\partial t}$ and integrating the resulting equation by parts.
Similarly, the energy associated with the thick structure problem is given by
\begin{align}
\nonumber 
&\displaystyle{ 
\frac{1}{2}\frac{d}{dt}\Big (\rho_s  \left\|\frac{\partial{\hat{\bf U}}}{\partial t} \right\|^2_{L^2(\hat\Omega_S)}
+2 \mu_s \|{\bf D}({\hat{\bf U}})\|^2_{L^2(\hat\Omega_S)}
+ \lambda_s \|\nabla\cdot{\hat{\bf U}}\|^2_{L^2(\hat\Omega_S)}}\\
&\qquad\qquad\qquad  + {\gamma} \|\hat{\bf U}\|^2_{L^2(\hat\Omega_S)}\Big )= \displaystyle{
-\int_0^L \boldsymbol \sigma^s {\bf e}_r\cdot\frac{\partial{\hat{\bf U}}}{\partial t},
}
\label{ThickEnergy}
\end{align}
which is obtained after multiplying the thick structure problem \eqref{linela} by structure velocity 
$\partial\hat{\bf U}/\partial t$ and integrating the resulting equation by parts.

To deal with the fluid equations we use the same approach: multiply the momentum equation \eqref{NS1}
by the fluid velocity and integrate by parts, using the incompressibility condition on the way, and the following identities:
\begin{eqnarray}
\int_{\Omega^f(t)} \frac{\partial \boldsymbol v}{\partial t} \boldsymbol v d\boldsymbol x &=& \frac{1}{2} \frac{d}{dt}\int_{\Omega^f(t)} |\boldsymbol v|^2 d\boldsymbol x - 
\frac{1}{2} \int_{\Gamma(t)} |\boldsymbol v|^2 \boldsymbol v \cdot \boldsymbol n^f dS, \label{ident1} 
\\
\int_{\Omega^f(t)} (\boldsymbol v \cdot \nabla) \boldsymbol v \cdot \boldsymbol v d\boldsymbol x &=& \frac{1}{2} \int_{\partial \Omega^f(t)} |\boldsymbol v|^2 \boldsymbol v \cdot \boldsymbol n^f dS, \label{ident2}
\end{eqnarray}
The first  identity is just the transport theorem, while 
the second one is obtained using integration by parts.
In \eqref{ident2}, the boundary integral over $\partial \Omega^f(t) =  \Gamma_b \cup \Gamma(t) \cup \Gamma_{in/out}$  is simplified as follows.
The portion corresponding to $\Gamma_b$ is zero due to the 
symmetry boundary condition, which implies ${\bf v}\cdot{\bf n}=0$ on $\Gamma_b$.
The portion corresponding to $\Gamma(t)$ is canceled with the same term appearing in the transport formula \eqref{ident1}. Finally,
the boundary terms on $\Gamma_{in/out}$ are absorbed by the dynamic pressure boundary conditions. 
Therefore, after multiplying~\eqref{NS1} with ${\bf v}$, integrating by parts and taking into account
the preceding discussion, we get:
\begin{eqnarray}
\nonumber
&\displaystyle{ \frac{\rho_f}{2}\frac{d}{dt}  ||\boldsymbol v||^2_{L^2(\Omega^f(t))}  + 2 \mu_f ||\boldsymbol D(\boldsymbol v)||^2_{L^2(\Omega(t))}
- \int_{\Gamma(t)} \boldsymbol\sigma \boldsymbol {n}^f \cdot \boldsymbol v\ dS}\\
&\displaystyle{ = - \int_0^R p_{in}(t) v_z|_{z=0} dr +   \int_0^R p_{out}(t) v_z|_{z=L} dr}.
\label{FluidEnergy}
\end{eqnarray}

The energy of the Koiter shell, given by \eqref{KoiterEnergy}, the energy of the thick structure, given by \eqref{ThickEnergy},
and the energy of the fluid, given by \eqref{FluidEnergy}, 
are combined into the total energy of the coupled problem by employing the dynamic and kinematic coupling conditions. 
To do this, first the fluid normal stress term in \eqref{FluidEnergy} is written in Lagrangian coordinates as follows:
\begin{equation}
\int_{\Gamma(t)} \boldsymbol\sigma \boldsymbol{n}^f \cdot \boldsymbol v\ dS
=\int_{\hat{\Gamma}}\widehat{\boldsymbol\sigma \boldsymbol{n}^f} \cdot \hat{\boldsymbol v} \ J \ d\hat{z},
\end{equation}
where $J$ is the Jacobian of the transformation from  Eulerian to Lagrangian coordinates.
From the kinematic and dynamic coupling conditions, the normal stress term is equal to:
\begin{equation}\label{weak_coupling}
\int_{\hat{\Gamma}} \widehat{\boldsymbol\sigma \boldsymbol {n}^f} \cdot \hat{\boldsymbol v}  \ J \ d\hat{z}
= -\int_{\hat{\Gamma}} \hat{\boldsymbol f} \cdot\frac{\partial \hat{\boldsymbol \eta}}{\partial t}  d\hat{z}-  \int_{\hat{\Gamma}} \boldsymbol \sigma^s \boldsymbol n^s  \cdot \frac{\partial \boldsymbol {U}}{\partial t} \bigg|_{\hat{\Gamma}}  d\hat{z}.
\end{equation}
This is used in the final step in which the three expressions for the energy of each separate physics sub-problem \eqref{KoiterEnergy},
\eqref{ThickEnergy}, and \eqref{FluidEnergy}, are added together, 
to obtain exactly the energy equality \eqref{FSIEnergy} of the coupled multi-physics FSI problem~\eqref{NS1} - \eqref{dynamic}
driven by the time-dependent dynamic pressure data.

Therefore, we have shown that if a smooth solution to the coupled fluid-structure interaction problem~\eqref{NS1} - \eqref{dynamic} exists, then it satisfies 
the  energy equality~\eqref{FSIEnergy}. This equality states that the rate of change of the kinetic energy of the fluid,  the kinetic energy of the multilayered structure,
and  the elastic energy of the structure, plus the viscous dissipation of the fluid,
is counter-balanced by the work done by the inlet and outlet data.

\subsection{The ALE framework}
To circumvent the difficulty associated with the fact that fluid domain changes in time we adopt the Arbitrary Lagranian Eulerian (ALE) method~\cite{hughes1981lagrangian,donea1983arbitrary,nobile2001numerical}. 
The ALE approach is based on introducing a family of arbitrary, smooth, homeomorphic mappings ${\cal{A}}_t$ defined on the reference domain $\hat{\Omega}^f$ such that, for each $t \in (t_0, T)$, ${\cal{A}}_t$ maps the reference domain $\hat{\Omega}^f$ into the current domain $\Omega^f(t)$:
$$\mathcal{A}_t : \hat{\Omega}^f  \rightarrow \Omega^f(t) \subset \mathbb{R}^2,  \quad \boldsymbol{x}=\mathcal{A}_t(\hat{\boldsymbol{x}}) \in \Omega^f(t), \quad \textrm{for} \; \hat{\boldsymbol{x}} \in \hat{\Omega}^f.$$

We will use the ALE mapping to deal with the deformation of the mesh, 
and to resolve the issues 
related to the approximation of the time-derivative
$\partial \boldsymbol v/\partial t \approx (\boldsymbol {v}(t^{n+1})-\boldsymbol {v}(t^{n}))/\Delta t,$ which
is not well-defined
due to the fact that $\Omega^f(t)$ depends on time. To rewrite the Navier-Stokes equations~\eqref{NS1} in the ALE framework
we employ the ALE mapping $\mathcal{A}_t$ to relate a function $f=f(\boldsymbol {x},t)$ defined on $\Omega^f(t) \times (0,T)$ 
to the corresponding function $\hat{f} := f \circ \mathcal{A}_t$ defined on $\hat{\Omega} \times (0,T)$
as
$$\hat{f}(\hat{\boldsymbol{x}},t) = f(\mathcal{A}_t(\hat{\boldsymbol{x}}),t).$$
Differentiating the latter equation with respect to time, we have
\begin{equation}
 \frac{\partial f}{\partial t}\bigg|_{\hat{\boldsymbol x}} =  \frac{\partial {f}}{\partial t}+\boldsymbol{w} \cdot \nabla f,
\end{equation}
where 
$\boldsymbol{w}=\displaystyle\frac{\partial \mathcal{A}_t(\hat{\boldsymbol x})}{\partial t}$ denotes the domain velocity. Note that $ \displaystyle\frac{\partial f}{\partial t}\bigg|_{\hat{\boldsymbol x}} $ denotes the 
time derivative of $f$ evaluated on the reference domain. 

Finally, system~\eqref{NS1}-\eqref{NS2} in the ALE framework reads as follows: 
Find $\boldsymbol{v}= (v_z,v_r)$ and $p$,  
with $\hat{\boldsymbol{v}}(\hat{\boldsymbol{x}},t) = \boldsymbol{v}(\mathcal{A}_t(\hat{\boldsymbol{x}}),t)$  such that
\begin{align}
&  \rho_f \bigg( \frac{\partial \boldsymbol{v}}{\partial t}\bigg|_{\hat{\boldsymbol x}}+ (\boldsymbol{{v}}-\boldsymbol{w}) \cdot \nabla \boldsymbol{v} \bigg)  = \nabla \cdot \boldsymbol\sigma^f(\boldsymbol v, p),  & \textrm{in}\; \Omega^f(t) \times (0,T), \\
& \nabla \cdot \boldsymbol{v} = 0  & \textrm{in}\; \Omega^f(t) \times (0,T),
\end{align}
with corresponding initial and boundary conditions.

\subsection{Weak formulation of the coupled FSI problem}
For all $t \in [0,T)$ define the following test function spaces  
\begin{eqnarray}
V^f(t) &=& \{\boldsymbol {\varphi}: \Omega^f(t) \rightarrow \mathbb{R}^2| \; \boldsymbol \varphi = \hat{\boldsymbol \varphi} \circ (\mathcal{A}_t)^{-1}, \hat{\boldsymbol \varphi} \in (H^1(\hat{\Omega}^f))^2,  \nonumber \\
& & \quad \quad  \varphi_r|_{r=0}=0, \; \boldsymbol \varphi|_{z=0,L}=0 \}, \label{Vf(t)} \\
Q(t) &=& \{q: \Omega^f(t) \rightarrow \mathbb{R}| \; q = \hat{q} \circ (\mathcal{A}_t)^{-1}, \hat{q} \in L^2(\hat{\Omega}^f) \}, \label{Q(t)} \\
\hat{V}^m &=& \{\hat{\boldsymbol {\zeta}}: (0,L) \rightarrow \mathbb{R}^2| \; \hat{\boldsymbol \zeta} \in (H_0^1(\hat{\Gamma}))^2 \}, \label{Vk} \\
\hat{V}^s &=& \{\hat{\boldsymbol {\psi}}: \hat{\Omega}^s \rightarrow \mathbb{R}^2| \; \hat{\boldsymbol \psi} \in (H^1(\hat{\Omega}^s))^2, \hat{\boldsymbol \psi} |_{z=0, L}=0, \hat\psi_z|_{\hat{\Gamma}^s_{ext}} = 0  \}, \label{Vs}
\end{eqnarray}
and introduce the function space 
\begin{equation}
V^{fsi} = \{(\boldsymbol \varphi, \hat{\boldsymbol \zeta}, \hat{\boldsymbol \psi}) \in V^f(t)\times \hat{V}^m \times \hat{V}^s| \; \boldsymbol \varphi|_{\Gamma(t)} = \boldsymbol \zeta, \; \hat{\boldsymbol \zeta} = \hat{\boldsymbol \psi}|_{\hat{\Gamma}} \},
\end{equation}
where $\boldsymbol \zeta := \hat{\boldsymbol \zeta} \circ (\mathcal{A}_t^{-1}|_{\Gamma(t)})$.
The variational formulation of the coupled fluid-structure interaction problem now reads: 
given $t \in (0, T)$, 
find 
$(\boldsymbol v,\hat{\boldsymbol \eta}, \hat{\boldsymbol U}, p) \in V^f(t)\times\hat{V}^m\times\hat{V}^s\times Q(t)$ 
such that~\eqref{kinematic} and~\eqref{kinematic1} hold and  such that for all $(\boldsymbol \varphi, \hat{\boldsymbol \zeta}, \hat{\boldsymbol \psi}, q) \in V^{fsi} \times Q(t)$
\begin{gather}
\rho_f \int_{\Omega^f(t)} \frac{\partial \boldsymbol v}{\partial t} \cdot \boldsymbol \varphi d\boldsymbol x+ \int_{\Omega^f(t)} (\boldsymbol v \cdot \nabla) \boldsymbol v \cdot \boldsymbol \varphi d \boldsymbol x  +2 \mu_f \int_{\Omega^f(t)} \boldsymbol D(\boldsymbol v) : \boldsymbol D(\boldsymbol \varphi) d \boldsymbol x 
\notag \\
- \int_{\Omega^f(t)} p \nabla \cdot \boldsymbol \varphi d \boldsymbol x+ \int_{\Omega^f(t)} q \nabla \cdot \boldsymbol v d \boldsymbol x+\rho_{m} h \int_0^L \frac{\partial^2 \hat\eta_z}{\partial t^2} \hat\zeta_z d\hat{z} + \rho_{m} h \int_0^L \frac{\partial^2 \hat\eta_r}{\partial t^2} \hat\zeta_r d\hat{z}
\notag \\
 -C_2 \int_0^L \frac{\partial \hat\eta_r}{\partial \hat{z}} \hat\zeta_z d\hat{z} +C_1 \int_0^L \frac{\partial \hat\eta_z}{\partial \hat{z}} \frac{\partial \hat\zeta_z}{\partial \hat{z}} d\hat{z}+C_0 \int_0^L \hat\eta_r \hat\zeta_r d\hat{z}  + C_2 \int_0^L  \frac{\partial \hat\eta_z}{\partial \hat{z}} \hat\zeta_r d\hat{z}
\notag \\
+\rho_{s} \int_{\hat{\Omega}^s} \frac{\partial^2 \hat{\boldsymbol U}}{\partial t^2} \cdot \boldsymbol \psi d\hat{\boldsymbol x}+2 \mu_{s} \int_{\hat{\Omega}^s} \boldsymbol D(\hat{\boldsymbol U}) : \boldsymbol D(\hat{\boldsymbol \psi}) d \hat{\boldsymbol x }+\lambda_{s} \int_{\hat{\Omega}^s} (\nabla \cdot \hat{\boldsymbol U}) (\nabla \cdot \hat{\boldsymbol \psi}) d \hat{\boldsymbol x}
\notag \\
 +\gamma \int_{\hat{\Omega}^s} \hat{\boldsymbol U} \cdot \hat{\boldsymbol \psi} d\hat{\boldsymbol x}= \int_0^R p_{in}(t) \varphi_z|_{z=0} dr
-  \int_0^R p_{out}(t) \varphi_z|_{z=L} dr,
\label{FSIweak}
\end{gather}
supplemented with initial conditions.

\subsection{Literature review}

The development of numerical solvers for fluid-structure interaction problems has become 
particularly active since the 1980's \cite{peskin1,peskin2,fauci1,fogelson,peskin3,peskin4,Griffith1,Griffith2,Griffith3,Griffith4,
donea1983arbitrary, hughes1981lagrangian,heil2004efficient,le2001fluid,leuprecht2002numerical,quarteroni2000computational,quaini2007semi,
baaijens2001fictitious,van2004combined,
fang2002lattice,feng2004immersed,krafczyk1998analysis,krafczyk2001two,
cottet2008eulerian,
figueroa2006coupled}.

\if 1 = 0
Among the most popular techniques are
the Immersed Boundary Method~\cite{peskin1,peskin2,fauci1,fogelson,peskin3,peskin4,Griffith1,Griffith2,Griffith3,Griffith4}
and the Arbitrary Lagrangian Eulerian (ALE) method
 \cite{donea1983arbitrary, hughes1981lagrangian,heil2004efficient,le2001fluid,leuprecht2002numerical,quarteroni2000computational,quaini2007semi}.
We further mention the Fictitious Domain Method in
combination with the mortar element method or ALE approach
\cite{baaijens2001fictitious,van2004combined}, and the methods recently proposed for the
use in the blood flow application such as
the Lattice Boltzmann method \cite{fang2002lattice,feng2004immersed,krafczyk1998analysis,krafczyk2001two},
the Level Set Method~\cite{cottet2008eulerian}
and the Coupled Momentum Method \cite{figueroa2006coupled}.
\fi

Until recently, only monolithic algorithms seemed  applicable to blood
flow simulations
\cite{figueroa2006coupled,gerbeau,nobile2008effective,zhao,bazilevs,bazilevs2}.
These algorithms are based on solving the entire nonlinear coupled problem as one monolithic system.
They are, however,
generally quite expensive in terms of computational time, programming time and
memory requirements, since they require solving a sequence of
strongly coupled problems using, e.g.,
the fixed point and Newton's methods
\cite{cervera,nobile2008effective,deparis,FernandezLifeV,heil2004efficient,matthies}.

The multi-physics nature of the blood flow problem
strongly suggest to employ partitioned (or staggered) numerical algorithms, where the coupled fluid-structure interaction
problem is separated into a fluid and a structure sub-problem.
The fluid and structure sub-problems are integrated in time in an alternating way, and the coupling conditions are enforced asynchronously. 
When the density of the structure is much larger than the density of the fluid, as is the case in aeroelasticity, 
it is sufficient to solve, at every time step, just one fluid sub-problem and one structure sub-problem to obtain a solution.
The classical {loosely-coupled} partitioned schemes of this kind typically use the structure velocity in the {fluid sub-problem} as
{Dirichlet} data for the fluid velocity (enforcing the no-slip boundary condition at the fluid-structure interface),
while in the { structure sub-problem} the structure is loaded by the fluid {normal stress} calculated in the fluid sub-problem.
These {{Dirichlet-Neumann}} loosely-coupled partitioned schemes work well for problems in which the structure is much heavier than the fluid. 
Unfortunately, when fluid and structure have comparable densities, which is  the case in the blood flow application, the simple strategy of separating the fluid from the structure suffers from severe stability 
issues \cite{causin2005added,Michler} associated with the added mass effect,
discussed in the introduction. 
The added mass effect reflects itself in Dirichlet-Neumann loosely coupled partitioned schemes
by having a poor approximation of the total energy of the coupled problem at every time step of the scheme. 
A partial solution to this problem is to iterate several times between the fluid and structure sub-solvers
at every time step until the energy of the continuous problem is well approximated.
These {strongly-coupled} partitioned schemes, however,  are computationally expensive and may suffer from
convergence issues for certain parameter values \cite{causin2005added}.

To get around these difficulties, and to retain the main advantages of loosely-coupled partitioned schemes
such as modularity, simple implementation, and low computational costs, several new loosely-coupled algorithms 
have been proposed recently. In general, they behave quite well for FSI problems containing a thin fluid-structure interface with mass
\cite{badia,Martina_paper1,MarSunLongitudinal,guidoboni2009stable,nobile2008effective,Fernandez1,Fernandez2,Fernandez3,Fernandez2006projection,astorino2009added,astorino2009robin,badia2008splitting,quaini2007semi,murea2009fast,deparis,deparis2}.

For FSI problem in which the structure is ``thick'' relative to the fluid, i.e., the thickness of the structure is comparable to the 
transverse dimension of the fluid domain, partitioned, loosely-coupled schemes
 are more difficult to construct. In fact, to the best of our knowledge, there have been no loosely-coupled, partitioned schemes
 proposed so far in literature for a class of FSI problems in hemodynamics that contain thick structure models to study the 
 elastodynamics of arterial walls. The closest works in this direction include a strongly-coupled partition scheme by Badia et al. in \cite{badia2009robin},
 an explicit scheme by Burman and Fern\'{a}ndez where certain ``defect-correction'' sub-iterations are necessary to achieve optimal accuracy
 \cite{burman2009stabilization}, and an operator-splitting scheme by Buka\v{c} et al. \cite{thick} where a version of Lie splitting was used to 
 design a modular, ``semi-partitioned'' scheme.

\if 1 = 0
 More precisely, in \cite{badia2009robin}, the authors construct a strongly-coupled partitioned scheme based on certain Robin-type coupling conditions. In addition to the classical Dirichlet-Neumann and Neumann-Dirichlet schemes, they
also propose a Robin-Neumann and a Robin-Robin scheme, that converge without relaxation, and
need a smaller number of sub-iteration between the fluid and the
structure in each time step than classical strongly-coupled schemes. 

In \cite{burman2009stabilization}, Burman and Fern\'{a}ndez propose an explicit scheme where 
the coupling between the fluid and a thick structure is enforced in a weak sense using Nitsche's approach \cite{hansbo2005nitsche}.
The formulation in \cite{burman2009stabilization} still suffers from stability issues related to the added mass effect, which were
corrected by adding a weakly consistent penalty term that includes pressure variations at the interface. The added term, however,
lowers the temporal accuracy of the scheme, which was then corrected by proposing a few defect-correction sub-iterations to achieve
optimal accuracy. 
\fi


Recently, a novel  loosely coupled partitioned scheme, called the Kinematically Coupled $\beta$-Scheme,
was introduced by Buka\v{c}, \v{C}ani\'{c} et al.~ in \cite{Martina_paper1,MarSunLongitudinal}, and
applied to FSI problems with thin elastic and viscoelastic structures, modeled by the membrane or shell equations.
This scheme successfully deals with stability problems
associated with the added mass effect
in a way different from those reported above. 
Stability is achieved by combining the structure inertia with the fluid sub-problem to mimic the energy balance 
of the continuous, coupled problem.
It was shown in \cite{SunBorMar} that the scheme is unconditionally stable even for the parameters associated with the blood flow applications.
Additionally, Muha and \v{C}ani\'{c} showed that a version of this scheme with $\beta=0$ converges to a weak solution of the fully nonlinear FSI problem \cite{BorSun}.
The case $\beta = 0$ considered in \cite{BorSun} corresponds to the classical kinematically-coupled scheme, first introduced in \cite{guidoboni2009stable}.
Parameter $\beta$  was introduced
in \cite{Martina_paper1} to increase the accuracy of the scheme. 
It was shown in \cite{Martina_paper1} that the accuracy of the kinematically-coupled $\beta$-scheme with $\beta = 1$
was comparable to that of monolithic scheme by Badia, Quaini, and Quarteroni in \cite{badia2008splitting}
when applied to the nonlinear benchmark FSI problem in hemodynamics, introduced by Formaggia et al. in \cite{formaggia2001coupling}.
A different approach to increasing the accuracy of the classical kinematically-coupled scheme was recently proposed
  by Fern{\'a}ndez et al. \cite{Fernandez1,Fernandez2,Fernandez3}. Their modified kinematically-coupled scheme, called
``the incremental displacement-correction scheme''  treats the structure displacement explicitly in the fluid sub-step 
and then corrects it in the structure sub-step. Fern{\'a}ndez et al. showed that the accuracy of the incremental displacement-correction scheme
is first-order in time. The results were obtained for a FSI problem involving a thin elastic structure.

These recent results indicate that the kinematically-coupled scheme and its modifications provide an appealing way to study 
multi-physics problems involving FSI.
Indeed, due to its simple implementation, modularity, and very good performance,
modifications of this scheme have been used by several authors to study different multi-physics problems involving FSI, such as
FSI between artery, blood flow, and cardiovascular device called stent \cite{BorSunStent},  
FSI with poroelastic structures \cite{Martina_Biot}, and FSI involving non-Newtonian fluids \cite{Lukacova,LukacovaCMAME}.

In the present paper we extend this scheme, in combination with the Arbitrary Lagrangian-Eulerian approach, 
 to study fluid structure interaction involving biological composite structures
such as arterial walls. Using a combination of analytical and numerical methods new results that reveal novel 
physical properties of the coupled FSI system are revealed. They include  regularizing effects of the fluid-structure interfere with mass,
and convergence, 
as the thickness of the thin fluid-structure interface with mass converges to zero,
to the FSI solution involving a single thick structure with the thickness
equal to the (fixed) combined thickness of the composite structure.

\section{The splitting algorithm}\label{splitting}
To numerically solve the fluid-multilayered structure interaction problem~\eqref{NS1}-\eqref{dynamic}
we propose a loosely coupled scheme based on the Lie operator splitting algorithm~\cite{glowinski2003finite}.
The proposed scheme is  an extension of the kinematically coupled $\beta-$scheme 
introduced in~\cite{Martina_paper1} to study the interaction between an incompressible, viscous fluid
and a thin structure modeled by the Koiter shell equations.
The scheme was shown to be unconditionally stable and therefore independent of the fluid and structure densities \cite{SunBorMar}.
The unconditional stability is achieved by including the thin structure inertia into the fluid sub-problem as Robin-type boundary condition. 

Following similar ideas, the proposed strategy for solving problem \eqref{NS1}-\eqref{dynamic} is to split the thin structure equations 
into the inertia part (plus the viscous part if the structure is viscoelastic), and the remaining elastodynamics part. The inertia part (plus
structural viscosity) is then used in the fluid sub-problem as Robin-type boundary condition, while the remaining
elastodynamics part is used in the structural sub-problem as a Robin-type boundary condition for the thick structure problem.
This strategy ensures unconditional stability of the scheme.
To increase the accuracy, the fluid stress will be split so that a portion corresponding to the fluid pressure is used as loading in the
structure sub-problem. 
To further increase the accuracy of the scheme, the fluid sub-problem can be split into a time-dependent Stokes problem
(with the Robin-type boundary condition mentioned above) plus the remaining advection sub-problem. This 
enables the use of conservative schemes to solve the advection sub-problem, if so desired.

\subsection{The Lie scheme and the first-order system}
To apply the Lie splitting scheme the problem must first be written as a first-order system in time:
\begin{eqnarray}\label{LieProblem}
   \frac{\partial \phi}{\partial t} + A(\phi) &=& 0, \quad \textrm{in} \ (0,T), \\
\phi(0) &=& \phi_0, 
\end{eqnarray}
where $A$ is an operator from a Hilbert space into itself. Operator $A$ is then split, in a non-trivial decomposition, as
\begin{equation}
 A = \sum\limits_{i=1}^I A_i.
\end{equation}
The Lie scheme consists of the following. Let $\Delta t>0$ be a time discretization step. Denote $t^n=n\Delta t$ and let $\phi^n$
be an approximation of $\phi(t^n).$ Set $\phi^0=\phi_0.$ Then, for $n \geq 0$ compute $\phi^{n+1}$ by solving
\begin{eqnarray}
   \frac{\partial \phi_i}{\partial t} + A_i(\phi_i) &=& 0 \quad \textrm{in} \; (t^n, t^{n+1}), \\
\phi_i(t^n) &=& \phi^{n+(i-1)/I}, 
\end{eqnarray}
and then set $\phi^{n+i/I} = \phi_i(t^{n+1}),$ for $i=1, \dots. I.$

This method is first-order accurate in time. 
More precisely, if~\eqref{LieProblem} is defined on a finite-dimensional space, and if the operators $A_i$ are smooth enough, 
then $\| \phi(t^n)-\phi^n \| = O(\Delta t)$~\cite{glowinski2003finite}.
In our case, operator $A$ that is associated with problem~\eqref{NS1}-\eqref{dynamic} 
will be split into a sum of three operators: 
\begin{enumerate}
 \item[ ]{\bf A1.} The thick structure elastodynamics problem with thin structure elastodynamics included as Robin-type boundary condition.
 \item[ ]{\bf A2(a).}  A time-dependent Stokes problem with thin structure inertia included as Robin boundary condition.  
  \item[ ]{\bf A2(b).} The fluid advection problem.
\end{enumerate}
Operator A1 defines a structure problem. Operators A2(a) and A2(b) define a fluid problem. 

To rewrite system~\eqref{NS1}-\eqref{dynamic} in first-order form, we introduce two new variables denoting the thin and thick structure velocities:
$$
\hat{\boldsymbol \xi}=\frac{\partial \hat{\boldsymbol \eta}}{\partial t} \ {\rm and}\ \hat{\boldsymbol V}=\frac{\partial \hat{\boldsymbol U}}{\partial t}.
$$
System~\eqref{NS1}-\eqref{dynamic} can be written in first-order ALE form as follows: 
Find $\boldsymbol{v}= (v_z,v_r)$, $p$, $\hat{\boldsymbol{\eta}}= (\hat\eta_z,\hat\eta_r)$, $\hat{\boldsymbol{\xi}}= (\hat\xi_z,\hat\xi_r),$ $\hat{\boldsymbol U} = (\hat{U}_z, \hat{U}_r)$ and $\hat{\boldsymbol V} = (\hat{V}_z, \hat{V}_r)$,
with $\hat{\boldsymbol{v}}(\hat{\boldsymbol{x}},t) = \boldsymbol{v}(\mathcal{A}_t(\hat{\boldsymbol{x}}),t),$  such that
\begin{subequations} \label{sys1}
\begin{align}
& \rho_f \bigg( \frac{\partial \boldsymbol{v}}{\partial t}\bigg|_{\hat{\boldsymbol x}}+ (\boldsymbol{{v}}-\boldsymbol{w}) \cdot \nabla \boldsymbol{v} \bigg)  = \nabla \cdot \boldsymbol\sigma^f(\boldsymbol v, p)   &  \textrm{in}\; \Omega^f(t) \times (0,T), \label{sys1a}\\
& \nabla \cdot \boldsymbol{v} = 0  & \textrm{in}\; \Omega^f(t) \times (0,T), \\
& \rho_{m} h \frac{\partial \hat\xi_z}{\partial t}-C_2 \frac{\partial \hat\eta_r}{\partial \hat{z}}-C_1 \frac{\partial^2 \hat\eta_z}{\partial\hat {z}^2}= \hat {f}_z  & \textrm{on} \; \hat{\Gamma} \times (0,T),  \label{fomem1}\\
& \rho_{m} h \frac{\partial \hat\xi_r}{\partial t}+C_0 \hat\eta_r +C_2 \frac{\partial \hat\eta_z}{\partial \hat{z}}   = \hat {f}_r & \textrm{on} \; \hat{\Gamma}\times (0,T),  \label{fomem2} \\
& \frac{\partial \hat{\boldsymbol\eta}}{\partial t} = \hat{\boldsymbol \xi}& \textrm{on} \; \hat{\Gamma}\times (0,T), \\
& \rho_{s} \frac{\partial \hat{\boldsymbol V}}{\partial t} + \gamma \hat{\boldsymbol U}= \nabla \cdot \boldsymbol \sigma^s( \hat{\boldsymbol U}) & \textrm{in} \; \hat{\Omega}^s \times (0,T), \\
& \frac{\partial \hat{\boldsymbol U}}{\partial t} = \hat{\boldsymbol V}   &  \textrm{in} \; \hat{\Omega}^s \times (0,T),
\end{align}
\end{subequations}
with the following coupling conditions at the fluid-structure interface:
\begin{eqnarray}
&  \boldsymbol \xi= \boldsymbol v|_{\Gamma(t)},  \quad \hat{\boldsymbol \eta} = \hat{\boldsymbol U}|_{\hat{\Gamma}}, \label{fokin}\\
 & \displaystyle{\rho_{m} h \frac{\partial (\widehat {v_z|_{\Gamma(t)}})}{\partial t}-C_2 \frac{\partial \hat\eta_r}{\partial \hat{z}}-C_1 \frac{\partial^2 \hat\eta_z}{\partial\hat {z}^2}+J\ \widehat{\boldsymbol{\sigma n}^f|_{\Gamma(t)} }\cdot \boldsymbol{e}_z + \boldsymbol \sigma^s\boldsymbol {n}^s|_{\hat{\Gamma}}\cdot \boldsymbol{e}_z} = 0, \\
& \displaystyle{\rho_{m} h \frac{\partial (\widehat {v_r|_{\Gamma(t)}})}{\partial t}+C_0 \hat\eta_r +C_2 \frac{\partial \hat\eta_z}{\partial \hat{z}}  +J\ \widehat{\boldsymbol{\sigma n}^f|_{\Gamma(t)} }\cdot \boldsymbol{e}_r+  \boldsymbol \sigma^s\boldsymbol {n}^s|_{\hat{\Gamma}}\cdot \boldsymbol{e}_r}  = 0.
\end{eqnarray}
This problem is supplemented with the boundary  and initial conditions presented in Section~\ref{description}.


\subsection{Details of the operator-splitting scheme}
The following splitting of the fluid stress $\widehat{\boldsymbol\sigma^f \boldsymbol {n}^f}$ will be used in the splitting algorithm:
$$
\widehat{{\boldsymbol\sigma}^f \boldsymbol{ n}^f} = \underbrace{\widehat{\boldsymbol\sigma^f \boldsymbol{n}^f}+\beta\widehat{p\boldsymbol{n}^f}}_{(Part \ I)}
\underbrace{- \beta \widehat{p\boldsymbol{n}^f}}_{(Part \ II)},
$$
where $\beta \in [0,1]$.
Part I of the fluid stress will be included in the Robin-type boundary condition
for the fluid sub-problem, while Part II of the fluid stress will be included in the Robin-type boundary condition for the thick structure equations. 
Details of the scheme are  as follows.

{\bf Problem A1:}  This problem involves solving the thick structure problem on $\hat\Omega_S$, 
with membrane elastodynamics serving as Robin-type boundary condition on $\hat\Gamma$. 
The Robin boundary condition also includes Part II of the normal fluid stress, evaluated at the previous time-step,
as part of the loading for the elastodynamics of the thin structure. 
After computing the thin and thick structure displacement, the domain velocity is computed and used in the 
fluid sub-problem in ALE mapping. 
The trace of the fluid velocity from the previous time step serves as initial data for the structure velocity in this problem.
Thus, the structure ``feels'' the fluid via the initial data and Part II of the normal fluid stress.

To be consistent with the notation in the Lie splitting scheme, we write Problem A1 in first-order form, i.e., as a first-order 
system in time, including all the variables in the problem. However, by using the kinematic coupling condition,
one can simplify the expression for the boundary condition on $\hat\Gamma$, as presented in Remark~1 below, 
in which case the present problem becomes a pure structure problem, driven by the fluid velocity initial data and fluid normal stress.

The problem reads:
Find $\boldsymbol{v},\hat{\boldsymbol \eta}, \hat{\boldsymbol \xi},\hat{\boldsymbol V},$ and $\hat{\boldsymbol U}$, with $\hat{\boldsymbol{v}}(\hat{\boldsymbol{x}},t) = \boldsymbol{v}(\mathcal{A}_t(\hat{\boldsymbol{x}}),t)$
such that:

\begin{equation}
\left.
\begin{array}{l}
\displaystyle{\rho_{s} \frac{\partial \hat{\boldsymbol V}}{\partial t} + \gamma \hat{\boldsymbol U}= \nabla \cdot \boldsymbol \sigma^s( \hat{\boldsymbol U})} \\
\displaystyle{\frac{\partial \hat{\boldsymbol U}}{\partial t} = \hat{\boldsymbol V}}
\end{array}
\right\} 
\ \textrm{in} \; \hat{\Omega}^s \times (t^n,t^{n+1}),
\end{equation}
satisfying the following conditions on $\hat\Gamma\times (t^n,t^{n+1})$:

\begin{align*}
 & \rho_{m} h \frac{\partial (\widehat {v_z|_{\Gamma(t)}})}{\partial t}-C_2 \frac{\partial \hat\eta_r}{\partial \hat{z}}-C_1 \frac{\partial^2 \hat\eta_z}{\partial\hat {z}^2}+ \boldsymbol \sigma^s\boldsymbol {n}^s|_{\hat{\Gamma}}\cdot \boldsymbol{e}_z 
= J^n \beta \widehat{p^{n} \boldsymbol n^f|_{\Gamma(t)}}\big)\cdot \boldsymbol{e}_z, \\
& \rho_{m} h \frac{\partial (\widehat {v_r|_{\Gamma(t)}})}{\partial t}+C_0 \hat\eta_r +C_2 \frac{\partial \hat\eta_z}{\partial \hat{z}} +  \boldsymbol \sigma^s\boldsymbol {n}^s|_{\hat{\Gamma}}\cdot \boldsymbol{e}_r  
= J^n \beta \widehat{p^{n} \boldsymbol n^f|_{\Gamma(t)}}\big)\cdot \boldsymbol{e}_r,\\
& \frac{\partial \hat{\boldsymbol\eta}}{\partial t} = \hat{\boldsymbol \xi}, \quad
\boldsymbol \xi= \boldsymbol v|_{\Gamma(t)},  \quad \hat{\boldsymbol \eta} = \hat{\boldsymbol U}|_{\hat{\Gamma}}, 
\end{align*}
and satisfying the following boundary conditions at $z=0,L$ and at $\hat\Gamma^s_{ext}$:
\begin{equation*}
\hat{\boldsymbol U}|_{z=0,L} = 0,
\end{equation*}
\begin{equation*}
 \hat{U}_z =  0, \quad \boldsymbol n^s_{ext} \cdot \boldsymbol \sigma^s\boldsymbol { n}^s_{ext} =  0 \quad \textrm{on} \; \hat{\Gamma}^s_{ext} \times (t^n,t^{n+1}).
\end{equation*}
The initial conditions are given by the values of the unknown functions at the end of the previous time step:
\begin{equation*}
\boldsymbol{v}(t^n)=\boldsymbol{v}^{n}, \; \hat{\boldsymbol \eta}(t^n) = \hat{\boldsymbol \eta}^{n}, \; \hat{\boldsymbol \xi}(t^n) = \hat{\boldsymbol \xi}^{n}, \; \hat{\boldsymbol{U}}(t^n)=\hat{\boldsymbol U}^{n}, \hat{\boldsymbol{V}}(t^n)=\hat{\boldsymbol V}^{n}.
\end{equation*}
The fluid velocity remains unchanged:
\begin{align*}
\frac{\partial \boldsymbol{v}}{\partial t}\bigg|_{\hat{\boldsymbol x}} = 0   \  \textrm{in} \; \Omega^f(t^n)\times(t^n, t^{n+1}). 
\end{align*}

The values of the just calculated variables are set to be $\boldsymbol{v}^{n+1/3}=\boldsymbol{v}(t^{n+1})$, $\hat{\boldsymbol \eta}^{n+1/3}=\hat{\boldsymbol \eta}(t^{n+1}), \; \hat{\boldsymbol \xi}^{n+1/3}=\hat{\boldsymbol \xi}(t^{n+1}), \; \hat{\boldsymbol U}^{n+1/3}=\hat{\boldsymbol U}(t^{n+1}),\; \hat{\boldsymbol V}^{n+1/3}=\hat{\boldsymbol V}(t^{n+1}), \; p^{n+1}=p(t^{n+1}).$ 

Using the just-calculated values of the thin structure displacement and velocity we  compute the  ALE velocity $\boldsymbol w$
by defining  $\mathcal{A}_{t^{n+1}}$ to be the harmonic extension of the structure displacement $\hat{\boldsymbol \eta}^{n+1/3}$ at the fluid-structure interface onto the whole domain $\hat\Omega^f$
\begin{eqnarray*}
 \Delta \mathcal{A}_{t^{n+1}} &=& 0 \quad \rm{in} \; \hat{\Omega}^f, \\
\mathcal{A}_{t^{n+1}} |_{\hat{\Gamma}} &=& \hat{\boldsymbol \eta}^{n+1/3}, \\
\mathcal{A}_{t^{n+1}} |_{\partial \hat{\Omega}^f\backslash \hat{\Gamma}}&=&0.
\end{eqnarray*}
Then, we define 
$
\boldsymbol w^{n+1} =\displaystyle \frac{\partial \mathcal{A}_{t^{n+1}} }{\partial t} = \displaystyle \frac{\partial \boldsymbol x}{\partial t}
\approx \frac{\boldsymbol x^{n+1} - \boldsymbol x^n}{\Delta t}
$, which remains unchanged in Problems 2(a) and 2(b).

\vskip 0.1in
\noindent
{\bf Remark 1.}
Note that using the kinematic coupling conditions, we can rewrite the membrane equations in the following way:
\begin{align*}
  & \rho_{m} h \frac{\partial \hat {V}_z}{\partial t}-C_2 \frac{\partial \hat{U}_r}{\partial \hat{z}}-C_1 \frac{\partial^2 \hat{U}_z}{\partial\hat {z}^2}+ \boldsymbol \sigma^s\boldsymbol {n}^s|_{\hat{\Gamma}}\cdot \boldsymbol{e}_z = J^n \beta \widehat{p^{n} \boldsymbol n^f|_{\Gamma(t)}}\big)\cdot \boldsymbol{e}_z, & \textrm{on} \; \hat{\Gamma}\times (t^n,t^{n+1}),\\
& \rho_{m} h \frac{\partial \hat {V}_r}{\partial t}+C_0 \hat{U}_r +C_2 \frac{\partial \hat{U}_z}{\partial \hat{z}} +  \boldsymbol \sigma^s\boldsymbol {n}^s|_{\hat{\Gamma}}\cdot \boldsymbol{e}_r  = J^n \beta \widehat{p^{n} \boldsymbol n^f|_{\Gamma(t)}}\big)\cdot \boldsymbol{e}_r & \textrm{on} \; \hat{\Gamma}\times (t^n,t^{n+1}).
\end{align*}
In this way the membrane equations become a Robin-type boundary condition for the thick structure problem, 
and Problem A1 becomes a pure structure problem. 

{\bf Problem A2(a).} This problem involves solving a time dependent Stokes problem with Robin-type boundary conditions involving the structure inertia
and Part I of the fluid stress. This problem is solved on a fixed fluid domain $\Omega^f(t^{n})$,
which is determined by the location of the structure position calculated in the previous time step.
To increase the accuracy, one can also use the just-updated structure position and solve the fluid sub-problem 
on domain $\Omega^f(t^{n+1})$. In the proof of stability of this scheme, the use of  $\Omega^f(t^{n})$ is more convenient. 

Similarly as before, we formally write the problem in terms of all the unknowns in the coupled problem. However,
after a close inspection of the problem, one can see that this is a pure fluid sub-problem, which communicates
with the structure problem via the thin structure inertia in the Robin-type boundary condition on $\hat\Gamma$, 
the updated fluid domain obtained from the structure problem in the previous time-step, and via the initial data 
for the trace of the fluid velocity on $\hat\Gamma$, given by the velocity of the thin structure, 
obtained from the just calculated structure step. 

 The problem reads as follows:
 Find $\boldsymbol{v}, p, \hat{\boldsymbol \eta}, \hat{\boldsymbol \xi},\hat{\boldsymbol V}$ and $\hat{\boldsymbol U}$, with $\hat{\boldsymbol{v}}(\hat{\boldsymbol{x}},t) = \boldsymbol{v}(\mathcal{A}_t(\hat{\boldsymbol{x}}),t)$
such that for $t\in (t^n, t^{n+1})$, with $p^n$ obtained at the previous time step, the following holds:

\begin{equation}
\left.
\begin{array}{rcl}
\displaystyle{
\rho_f \frac{\partial \boldsymbol{v}}{\partial t}\bigg|_{\hat{\boldsymbol x}}} &=&\nabla \cdot \boldsymbol{\sigma}^f  \\
\displaystyle{
\nabla \cdot \boldsymbol{v}}&=&0
\end{array}
\right\} \  \textrm{in} \; \Omega^f(t^{n}) \times(t^n, t^{n+1}),
\end{equation}
with the Robin-type boundary condition on $\hat\Gamma \times (t^n,t^{n+1})$ given by:
\begin{equation*}\label{step1}
\rho_{m} h \frac{\partial (\widehat {\boldsymbol{v}|_{\Gamma(t)}})}{\partial t}+J\big(\ \widehat{\boldsymbol{\sigma}^f \boldsymbol{n}^f|_{\Gamma(t)} }+\beta \widehat{p^n \boldsymbol n^f|_{\Gamma(t)}}\big) = 0.
\end{equation*}
This is supplemented with the following inlet and outlet boundary conditions on $(t^n,t^{n+1})$:
\begin{equation*}
  \boldsymbol{v}(0,R,t) = \boldsymbol{v}(L,R,t) = 0, 
\end{equation*}
\begin{equation*}
\boldsymbol{\sigma}^f \boldsymbol {n}^f_{in} = -p_{in}(t)\boldsymbol{n}^f_{in}\  {\rm on}\ \Gamma^f_{\rm in}, \; \; \boldsymbol{\sigma}^f \boldsymbol {n}^f_{out} = -p_{out}(t)\boldsymbol{n}^f_{out}  \ {\rm on}\  \Gamma^f_{\rm out},
\end{equation*}
and the following symmetry boundary conditions at $r = 0$:
\begin{equation*}
  \frac{\partial v_z}{\partial r}(z,0,t) = \quad v_r(z,0,t) = 0 \quad \textrm{on} \; (0,L)\times(t^n,t^{n+1}).
\end{equation*}
The location of the structure remains unchanged:
\begin{align*}
& \frac{\partial \hat{\boldsymbol\eta}}{\partial t} = 0 & \textrm{on} \; \hat{\Gamma}\times (t^n,t^{n+1}), \\
& \frac{\partial \hat{\boldsymbol V}}{\partial t} = \frac{\partial \hat{\boldsymbol U}}{\partial t}= 0 &  \textrm{in} \; \hat{\Omega}^s\times(t^n, t^{n+1}). 
\end{align*}
The initial conditions are given by:\\
$
\boldsymbol{v}(t^n)=\boldsymbol{v}^{n+1/3}, \; \hat{\boldsymbol \eta}(t^n) = \hat{\boldsymbol \eta}^{n+1/3}, \; \hat{\boldsymbol \xi}(t^n) = \hat{\boldsymbol \xi}^{n+1/3}, \; \hat{\boldsymbol{U}}(t^n)=\hat{\boldsymbol U}^{n+1/3}, 
$ and
$\hat{\boldsymbol{V}}(t^n)=\hat{\boldsymbol V}^{n+1/3}.$

After a solution to this problem is calculated, set $\boldsymbol{v}^{n+2/3}=\boldsymbol{v}(t^{n+1}), \; \hat{\boldsymbol \eta}^{n+2/3}=\hat{\boldsymbol \eta}(t^{n+1}), \; \hat{\boldsymbol \xi}^{n+2/3}=\hat{\boldsymbol \xi}(t^{n+1}), \; \hat{\boldsymbol U}^{n+2/3}=\hat{\boldsymbol U}(t^{n+1}),\; \hat{\boldsymbol V}^{n+2/3}=\hat{\boldsymbol V}(t^{n+1}), \; p^{n+1}=p(t^{n+1}).$

{\bf Problem A2(b).} Solve the fluid and ALE advection sub-problem defined on the fixed domain $\Omega^f(t^n)$,
with the domain velocity $\boldsymbol w^{n+1}$ just calculated in Problem A1. 

The problem reads:
Find $\boldsymbol{v},\hat{\boldsymbol \eta}, \hat{\boldsymbol \xi},\hat{\boldsymbol V}$ and $\hat{\boldsymbol U}$, with $\hat{\boldsymbol{v}}(\hat{\boldsymbol{x}},t) = \boldsymbol{v}(\mathcal{A}_t(\hat{\boldsymbol{x}}),t)$ and $\boldsymbol w^{n+1}$ computed in Problem A1,
such that:
\begin{equation*}
\frac{\partial \boldsymbol{v}}{\partial t}\bigg|_{\hat{\boldsymbol x}} + (\boldsymbol{v}^{n}-\boldsymbol{w}^{n+1}) \cdot \nabla \boldsymbol{v}= 0
\quad  \textrm{in} \; \Omega^f(t^{n})\times (t^n, t^{n+1}), 
\end{equation*}
with the lateral boundary condition on $\hat{\Gamma}\times (t^n, t^{n+1})$ given by
\begin{equation*}
 \rho_{m} h \frac{\partial (\widehat {\boldsymbol{v}|_{\Gamma(t)}})}{\partial t} = 0,
 \end{equation*}
and with the inlet/outlet boundary conditions defined on the portion of the boundary for which the slope of the characteristics associated 
with the advection problem, defined by the sign of $(\boldsymbol{v}^{n}-\boldsymbol{w}^{n+1})\cdot \boldsymbol{n}$, is pointing
inward, toward the fluid domain : 
$$\boldsymbol{v}=\boldsymbol{v}^{n+2/3}, \  \; \textrm{on} \; \Gamma_{-}^{n+2/3}, \; \textrm{where}$$ 
$$\Gamma_{-}^{n+2/3} = \{\boldsymbol{x} \in \mathbb{R}^2 | \boldsymbol{x} \in \partial \Omega^f(t^{n}), (\boldsymbol{v}^{n}-\boldsymbol{w}^{n+1})\cdot \boldsymbol{n} <0 \}.$$ 
The structure does not change in this step:
\begin{align*}
& \frac{\partial \hat{\boldsymbol\eta}}{\partial t} = 0 & \textrm{on} \; \hat{\Gamma}\times (t^n,t^{n+1}), \\
& \frac{\partial \hat{\boldsymbol V}}{\partial t} =\frac{\partial \hat{\boldsymbol U}}{\partial t}=  0 &  \textrm{in} \; \hat{\Omega}^s\times(t^n, t^{n+1}), 
\end{align*}
The initial conditions are given by the solution obtained in Problem A2(a):\\
$
\boldsymbol{v}(t^n)=\boldsymbol{v}^{n+2/3}, \; \hat{\boldsymbol \eta}(t^n) = \hat{\boldsymbol \eta}^{n+2/3}, \; \hat{\boldsymbol \xi}(t^n) = \hat{\boldsymbol \xi}^{n+2/3}, \; \hat{\boldsymbol{U}}(t^n)=\hat{\boldsymbol U}^{n+2/3}, 
$ 
and
$\hat{\boldsymbol{V}}(t^n)=\hat{\boldsymbol V}^{n+2/3}.$

Then set $\boldsymbol{v}^{n+1}=\boldsymbol{v}(t^{n+1}), \; \hat{\boldsymbol \eta}^{n+1}=\hat{\boldsymbol \eta}(t^{n+1}), \; \hat{\boldsymbol \xi}^{n+1}=\hat{\boldsymbol \xi}(t^{n+1}), \; \hat{\boldsymbol U}^{n+1}=\hat{\boldsymbol U}(t^{n+1}),\; \hat{\boldsymbol V}^{n+1}=\hat{\boldsymbol V}(t^{n+1}).$ \\

After this step we update the fluid domain position via
 \begin{equation*}
 \Omega^f(t^{n+1}) = \mathcal{A}_{t^{n+1}}(\hat{\Omega}^f) .
\end{equation*}
Then set  $t^n=t^{n+1}$ and return to Problem A1.

\section{Discretized scheme in weak form}\label{discretization}

As mentioned above, discretization in time 
is obtained by subdividing the time interval $(0,T)$ into $N$ sub-intervals of width $\Delta t$, 
and denoting $t^n = n \Delta t$, where $n \le N$. The Backward Euler scheme is used to approximate the time derivatives on each $(t^{n+1}, t^n)$. 
For the space discretization, we use the finite element method approach, based on a conforming FEM triangulation with maximum triangle diameter $k$. 
For this purpose, we define the finite element spaces 
$V^f_k(t^n) \subset V^f(t^n), Q_k(t^n) \subset Q(t^n), \hat{V}^m_k \subset \hat{V}^m$ and $ \hat{V}^s_k \subset \hat{V}^s$,
and introduce the following bilinear forms
\begin{eqnarray*}
 a^n_f(\boldsymbol v, \boldsymbol \varphi) &:=& 2 \mu_f \int_{\Omega^f(t^n)} \boldsymbol D(\boldsymbol v) : \boldsymbol D(\boldsymbol \varphi) d \boldsymbol x,   \label{af}\\
 b^n_f(p, \boldsymbol \varphi) &:=&  \int_{\Omega^f(t^n)} p \nabla \cdot \boldsymbol \varphi d \boldsymbol x, \label{bf} \\
 \hat{a}_m (\hat{\boldsymbol\eta}, \hat{\boldsymbol\zeta}) &:=&  C_1 \int_0^L \frac{\partial \hat\eta_z}{\partial \hat{z}} \frac{\partial \hat\zeta_z}{\partial \hat{z}} d\hat{z}+C_0 \int_0^L \hat\eta_r \hat\zeta_r d\hat{z}   \\
 &  + & C_2 \int_0^L  \frac{\partial \hat\eta_z}{\partial \hat{z}} \hat\zeta_r d\hat{z} 
 -C_2 \int_0^L \frac{\partial \hat\eta_r}{\partial \hat{z}} \hat\zeta_z d\hat{z} , \\
  \hat{a}_e(\hat{\boldsymbol U}, \hat{\boldsymbol \psi}) &:=& 2 \mu_s \int_{\hat{\Omega}^s} \boldsymbol D(\hat{\boldsymbol U}) : \boldsymbol D(\hat{\boldsymbol \psi}) d \hat{\boldsymbol x} + \lambda_s \int_{\hat{\Omega}^s} (\nabla \cdot \hat{\boldsymbol U})(\nabla \cdot \hat{\boldsymbol \psi})  d \hat{\boldsymbol x}   \\
  & + & \gamma \int_{\hat{\Omega}^s} \hat{\boldsymbol U} \cdot \hat{\boldsymbol \psi}  d \hat{\boldsymbol x}. 
  \label{ae}
 \end{eqnarray*}
Discretization of the operator splitting scheme discussed above, in weak form, is given by the following:
\vskip 0.1in
 \textbf{Problem A1. (The structure problem)} To discretize the structure problem in time we use the second order Newmark scheme. The problem reads as follows: Find $(\hat{\boldsymbol U}_k^{n+1/3}, \hat{\boldsymbol V}_k^{n+1/3}) \in \hat{V}_h^s\times \hat{V}_k^s$  such that for all $(\hat{\boldsymbol \psi}_k,\hat{\boldsymbol \phi}_k) \in \hat{V}_k^s \times \hat{V}_k^s$ 
\begin{gather}
\displaystyle{\rho_{s} \int_{\hat{\Omega}^s} \frac{\hat{\boldsymbol V}_k^{n+1/3}-\hat{\boldsymbol V}_k^{n}}{\Delta t}\cdot \hat{\boldsymbol \psi}_k d \hat{\boldsymbol x}+\gamma \int_{\hat{\Omega}^s} \frac{\hat{\boldsymbol U}_k^{n}+\hat{\boldsymbol U}_k^{n+1/3}}{2}\cdot \hat{\boldsymbol \psi}_k d\hat{\boldsymbol x}} 
\notag \\
 +\rho_m h \int_0^L \frac{\hat{\boldsymbol V}_k^{n+1/3}-\hat{\boldsymbol V}_{k}^n}{\Delta t} \cdot \hat{\boldsymbol\psi}_{k}d \hat{z}
 +\hat{a}_m(\frac{\hat{\boldsymbol U}^{n}_{k}+\hat{\boldsymbol U}^{n+1/3}_{k}}{2},\hat{\boldsymbol \psi}_{k})
\notag \\
+\hat{a}_e(\frac{\hat{\boldsymbol U}_k^{n}+\hat{\boldsymbol U}_k^{n+1/3}}{2}, \hat{\boldsymbol \psi}_k)+\displaystyle{ \rho_{s} \int_{\hat{\Omega}^s} (\frac{\hat{\boldsymbol V}_k^{n}+\hat{\boldsymbol V}_k^{n+1/3}}{2}- \frac{\hat{\boldsymbol U}_k^{n+1/3}-\hat{\boldsymbol U}_k^{n}}{\Delta t} ) \cdot \hat{\boldsymbol \phi}_k d \hat{\boldsymbol x}} \notag
\end{gather}
\begin{gather}
+\displaystyle{ \rho_{m}h \int_0^L (\frac{\hat{\boldsymbol V}_k^{n}+\hat{\boldsymbol V}_k^{n+1/3}}{2}- \frac{\hat{\boldsymbol U}_k^{n+1/3}-\hat{\boldsymbol U}_k^{n}}{\Delta t} ) \cdot \hat{\boldsymbol\phi}_k d \hat{ z}} \notag \\
 = \int_0^L J^n \beta \widehat{p^{n} \boldsymbol n^f|_{\Gamma(t)}}\cdot \hat{\boldsymbol \psi}_k d \hat{z}. \label{step1dfull}
\end{gather}
 Note that in this step we take all the kinematic coupling conditions into account: 
 initially we set $\hat{\boldsymbol V}_k^{n}|_{\Gamma}=\hat{\boldsymbol\xi}_k^{n}=\widehat{\boldsymbol v_k^{n}|_{\Gamma(t^n)}}$;
 then, once $\hat{\boldsymbol U}_k^{n+1/3}$ and $\hat{\boldsymbol V}_k^{n+1/3}$ are computed, 
 $\hat{\boldsymbol\eta}_k^{n+1/3}, \hat{\boldsymbol\xi}_k^{n+1/3}$ and $\widehat{\boldsymbol v_k^{n+1/3}|_{\Gamma(t^n)}}$ 
 can be recovered via $\hat{\boldsymbol \eta}^{n+1/3}_k = \hat{\boldsymbol U}_k^{n+1/3}|_{\Gamma}$ and $\hat{\boldsymbol \xi}^{n+1/3}_k = \widehat{\boldsymbol v^{n+1/3}_k|_{\Gamma(t^n)}}=\hat{\boldsymbol V}_k^{n+1/3}|_{\Gamma}$. 
 
 In this step ${\partial {\boldsymbol v}_k}/{{\partial t}} = 0,$ and so $\boldsymbol v_k^{n+1/3}=\boldsymbol v_k^{n}.$
 
At the end of this step we compute the new fluid domain velocity by finding $\mathcal{A}_{t^{n+1}} \in (H^1(\hat{\Omega}^f))^2$ with $\mathcal{A}_{t^{n+1}}|_{\hat{\Gamma}} = \hat{\boldsymbol \eta}^{n+1/3}$ and $\mathcal{A}_{t^{n+1}}|_{\partial \hat{\Omega}^f \backslash\hat{\Gamma}} = 0$ , such that
\begin{gather}
\int_{\hat{\Omega}^f} \nabla \mathcal{A}_{t^{n+1}} \cdot \nabla \boldsymbol \zeta d \hat{\boldsymbol x} =  0, \quad  \forall \boldsymbol \zeta \in (H^1_0(\hat{\Omega}^f))^2.
\end{gather}
Then, the new fluid domain velocity is given by $$
\boldsymbol w^{n+1} =\displaystyle \frac{\partial \mathcal{A}_{t^{n+1}} }{\partial t} 
\approx \frac{\mathcal{A}_{t^{n+1}}(\hat{\boldsymbol x}) - \boldsymbol x^n}{\Delta t}.
$$

 \vskip 0.1in
  \textbf{Problem A2(a). (The time dependent Stokes problem)}  We discretize Problem A2(a) using the Backward Euler scheme,  giving rise to the following weak formulation: Find $(\boldsymbol{v}_k^{n+2/3},p_k^{n+1}) \in V^{f}_k(t^{n}) \times Q_k(t^{n})$ such that for all $(\boldsymbol \varphi_k,q_k) \in V_k^{f}(t^{n}) \times Q_k(t^{n})$ 
\begin{gather}
\rho_f \int_{\Omega^f(t^{n})} \frac{\boldsymbol v_k^{n+2/3}-\boldsymbol v_k^{n+1/3}}{\Delta t} \cdot \boldsymbol \varphi_k d\boldsymbol x
 + a^{n}_f(\boldsymbol v_k^{n+2/3}, \boldsymbol \varphi_k)- b^{n}_f(p_k^{n+1}, \boldsymbol \varphi_k)
\notag \\
+  b^{n}_f(q_k, \boldsymbol v_k^{n+2/3})  +\rho_m h \int_{0}^L \frac{\widehat{\boldsymbol v_k^{n+2/3}|_{\Gamma(t^n)}}-\widehat{\boldsymbol v_k^{n+1/3}|_{\Gamma(t^n)}}}{\Delta t} \cdot \hat{\boldsymbol \varphi}_k d \hat{x}
\notag \\
 = \int_0^R p_{in}(t^{n+1}) \varphi_{z,k}|_{z=0} dr - \int_0^R p_{out}(t^{n+1}) \varphi_{z,k}|_{z=L} dr \notag \\
-\int_{0}^L J^n \beta \widehat{p^{n} \boldsymbol n^f|_{\Gamma(t)}}\cdot \hat{\boldsymbol \varphi}_k d \hat{x}. \label{step2adfull}
 \end{gather}
  
 In this step the structure is not moving, and so $\hat{\boldsymbol U}_k^{n+2/3}=\hat{\boldsymbol U}_k^{n+1/3}$,  $\hat{\boldsymbol V}_k^{n+2/3}=\hat{\boldsymbol V}_k^{n+1/3}$,
$\hat{\boldsymbol \eta}^{n+2/3}_k=\hat{\boldsymbol \eta}^{n+1/3}_k$, and $\hat{\boldsymbol \xi}^{n+2/3}_k=\hat{\boldsymbol \xi}^{n+1/3}_k$.

 \vskip 0.1in
  \textbf{Problem A2(b). (The fluid and ALE advection)} 
 As in the previous step, we discretize Problem A2(b) using the Backward Euler scheme. The weak formulation reads as follows: Find $\boldsymbol{v}_k^{n+1} \in V^{f}_k(t^{n})$ with $\boldsymbol{v}_k^{n+1} = \boldsymbol{v}_k^{n+2/3}$ on $\Gamma_{-}^{n+2/3},$ such that for all $\boldsymbol \varphi_k \in \{\boldsymbol \varphi_k \in V_k^{f}(t^n) | \ \boldsymbol\varphi_k=0 \; \textrm{on} \; \Gamma_{-}^{n+2/3} \}$ 
\begin{gather}
\rho_f \int_{\Omega^f(t^{n})} \frac{\boldsymbol v_k^{n+1}-\boldsymbol v_k^{n+2/3}}{\Delta t} \cdot \boldsymbol \varphi_k d\boldsymbol x
+\rho_f \int_{\Omega^f(t^n)}(({\boldsymbol v}_k^n-{\boldsymbol w}_k^{n+1})\cdot\nabla){\boldsymbol v}_k^{n+1}\cdot \boldsymbol {\varphi}_k =0.
\label{step2bdfull}
 \end{gather}
 Note that here we use ${\boldsymbol w}_k^{n+1}$ computed at the very end of  Problem {\bf A1}.

 In this step the structure is not moving, and so $$\hat{\boldsymbol U}_k^{n+1}=\hat{\boldsymbol U}_k^{n+2/3}, \hat{\boldsymbol V}_k^{n+1}=\hat{\boldsymbol V}_k^{n+2/3},
\hat{\boldsymbol\eta}^{n+1}_k=\hat{\boldsymbol\eta}^{n+2/3}_k,  \hat{\boldsymbol\xi}^{n+1}_k=\hat{\boldsymbol\xi}^{n+2/3}_k.$$

 Finally, we update the fluid domain position via
 \begin{equation*}
 \Omega^f(t^{n+1}) = (\boldsymbol I+\mathcal{A}_{t^{n+1}})(\hat{\Omega}^f) .
\end{equation*}
Then, set  $t^n=t^{n+1}$ and return to Problem {\bf A1}.

\section{Unconditional stability of the scheme for $\beta = 0$}\label{stability}
 
We show an energy estimate that is associated with
unconditional stability for the scheme
 when $\beta = 0$. 
We show that the sum of the total energy of the coupled problem plus 
viscous dissipation of the coupled discretized problem
is bounded by the discrete energy of initial data and the work done by the inlet and outlet dynamic pressure data.
In contrast with similar results appearing in literature which consider simplified models without fluid advection, 
and/or linearized fluid-structure coupling calculated at a fixed fluid domain boundary
\cite{causin2005added,SunBorMar,badia,Fernandez3,burman2009stabilization}, 
in this manuscript we prove the desired energy estimate for the full, nonlinear FSI problem, 
that includes fluid advection and nonlinear coupling at the moving fluid-structure interface.

To simplify analysis, the following assumptions that do not influence stability of the scheme will be considered:
\begin{description}
\item[1.] Only radial displacement of the fluid-structure interface is allowed, i.e., $\eta_z = U_z|_{\hat\Gamma} = 0$. In that case, we consider the following equation for the radial membrane displacement
 $$\rho_{m} h \frac{\partial^2 \hat\eta_r}{\partial t^2}+C_0 \hat\eta_r    = \hat {f}_r. $$
 The FSI problem with this boundary condition is well-defined. 
\item[2.] The problem is driven by the dynamic inlet and outlet pressure data, and the flow enters and leaves
the fluid domain parallel to the horizontal axis: 
$$
p+\frac{\rho_f}{2} |\boldsymbol v|^2 = p_{in/out}(t), \  v_r = 0, \ {\rm on}\ \Gamma^f_{in/out}.
$$
\end{description}
These assumptions do not influence stability issues related to the added mass effect
associated with partitioned schemes in FSI in blood flow.

To simplify matters, instead of splitting the Navier-Stokes equations into the time dependent Stokes problem and the advection problem, we keep the fluid problem uncoupled in the stability analysis. Therefore, our splitting scheme becomes:
\begin{itemize}
	\item[{\bf A1.}] The thick structure problem with membrane elastodynamics as Robin-type boundary condition.
	\item[{\bf A2.}] The fluid problem modeled by the Navier-Stokes equations in ALE form, 
	                              with structure inertia included as Robin boundary condition for the fluid problem.
	\end{itemize}

We consider the case $\beta = 0$ here.
Numerical results presented in \cite{Martina_paper1} indicate that only accuracy is affected by changing
the parameter $\beta$ between $0$ and $1$, and not stability.
We mention a related work in \cite{SunBorMar} where unconditional stability of the kinematically-coupled $\beta$-scheme
for $\beta\in[0,1]$ was proved for a simplified, linearized FSI problem with a thin structure, using different techniques from those
presented here.

\subsection{Discretized scheme in weak form assuming only radial displacement of the fluid-structure interface}

There are only a few differences with the discretized scheme  presented above. The bilinear form $\hat{a}_m$ associated 
with the thin structure problem is now simpler and is given by:
\begin{equation*}
 \hat{a}_m (\hat{\eta}_r, \hat{\zeta}_r) = C_0 \int_0^L \hat\eta_r \hat\zeta_r d\hat{z}.
 \end{equation*}
 Furthermore, we consider the problem with $\beta = 0$, and the splitting into two sub-problems: the structure and the fluid sub-problem,
 where, unlike before, the fluid dissipation and advection are treated together in the same step. 
 Finally, since we can easily explicitly calculate an ALE mapping in this case, the calculations of the ALE mapping,
 ALE velocity and the Jacobian of the ALE mapping are presented explicitly at the end of this sub-section.
 They are used in the derivation of the energy estimate, presented in the next sub-section.
 
 For completeness, we present the two steps of the discretized scheme next.
\vskip 0.1in
 \textbf{Problem A1 (Structure). }  Find $(\hat{\boldsymbol U}_k^{n+1/2}, \hat{\boldsymbol V}_k^{n+1/2}) \in \hat{V}_k^s\times \hat{V}_k^s$  such that for all $(\hat{\boldsymbol \psi}_k,\hat{\boldsymbol \phi}_k) \in \hat{V}_k^s \times \hat{V}_k^s$ 
\begin{equation*}\label{S1discrete}
\displaystyle{\rho_{s} \int_{\hat{\Omega}^s} \frac{\hat{\boldsymbol V}_k^{n+1/2}-\hat{\boldsymbol V}_k^{n}}{\Delta t}\cdot \hat{\boldsymbol \psi}_k d \hat{\boldsymbol x}+\gamma \int_{\hat{\Omega}^s} \frac{\hat{\boldsymbol U}_k^{n}+\hat{\boldsymbol U}_k^{n+1/2}}{2}\cdot \hat{\boldsymbol \psi}_k d\hat{\boldsymbol x}} 
\end{equation*}
\begin{equation*}
 +\rho_m h \int_{\Gamma} \frac{\hat{V}_{r,k}^{n+1/2}-\hat{V}_{r,k}^n}{\Delta t} \hat{\psi}_{r,k}
 +\hat{a}_m(\frac{\hat{U}^{n}_{r,k}+\hat{U}^{n+1/2}_{r,k}}{2},\hat{\psi}_{r,k})
+\hat{a}_e(\frac{\hat{\boldsymbol U}_k^{n}+\hat{\boldsymbol U}_k^{n+1/2}}{2}, \hat{\boldsymbol \psi}_k)
\end{equation*}
\begin{equation*}
+\displaystyle{ \rho_{s} \int_{\hat{\Omega}^s} (\frac{\hat{\boldsymbol V}_k^{n}+\hat{\boldsymbol V}_k^{n+1/2}}{2}- \frac{\hat{\boldsymbol U}_k^{n+1/2}-\hat{\boldsymbol U}_k^{n}}{\Delta t} ) \cdot \hat{\boldsymbol \phi}_k d \hat{\boldsymbol x}} 
\end{equation*}
\begin{equation}
+\displaystyle{ \rho_{m}h \int_{\Gamma} (\frac{\hat{V}_{r,k}^{n}+\hat{V}_{r,k}^{n+1/2}}{2}- \frac{\hat{U}_{r,k}^{n+1/2}-\hat{U}_{r,k}^{n}}{\Delta t} ) \cdot \hat \phi_{r,k} d \hat{ z}} 
 = 0, \label{step1d}
\end{equation}
with the initial data such that
 $\hat{V}_{r,k}^{n}|_{\Gamma}=\hat{\xi}_{r,k}^{n}=\widehat{v_{r,k}^{n}|_{\Gamma(t^n)}}.$
  
 Then set $\hat{\eta}^{n+1/2}_{r,k} = \hat{U}_{r,k}^{n+1/2}|_{\Gamma}$ and $\hat{\xi}^{n+1/2}_{r,k} = \widehat{v^{n+1/2}_{r,k}|_{\Gamma(t^n)}}=\hat{V}_{r,k}^{n+1/2}|_{\Gamma}$, and, since the fluid velocity does not change,  $\boldsymbol v_k^{n+1/2}=\boldsymbol v_k^{n}.$

 \vskip 0.1in
  \textbf{Problem A2 (Fluid).} Find $(\boldsymbol{v}_k^{n+1},p_k^{n+1}) \in V^{f}_k(t^{n}) \times Q_k(t^{n})$ such that for all $(\boldsymbol \varphi_k,q_k) \in V_k^{f}(t^{n+1}) \times Q_k(t^{n})$:
\begin{equation*}
\rho_f \int_{\Omega^f(t^{n})} \frac{\boldsymbol v_k^{n+1}-\boldsymbol v_k^{n+1/2}}{\Delta t} \cdot \boldsymbol \varphi_k d\boldsymbol x+\frac {\rho_f}{2} \int_{\Omega^f(t^n)}(\nabla\cdot{\boldsymbol w}_k^{n+1}){\boldsymbol v}_k^{n+1}\cdot\boldsymbol \varphi_k d \boldsymbol x
 \end{equation*}
 \begin{equation*}
+\frac{\rho_f}{2}\int_{\Omega^f(t^n)}\left (({\boldsymbol v}_k^n-{\boldsymbol w}_k^{n+1})\cdot\nabla){\boldsymbol v}_k^{n+1}\cdot \boldsymbol {\varphi}_k
-(({\boldsymbol v}_k^n-{\boldsymbol w}_k^{n+1})\cdot\nabla)\boldsymbol{\varphi}_k\cdot{\boldsymbol v}_k^{n+1}\right )
\end{equation*}
 \begin{equation*}
 + a^{n}_f(\boldsymbol v_k^{n+1}, \boldsymbol \varphi_k)- b^{n}_f(p_k^{n+1}, \boldsymbol \varphi_k)
+  b^{n}_f(q_k, \boldsymbol v_k^{n+1}) 
 \end{equation*}
 \begin{equation*}
 +\rho_m h \int_{0}^L \frac{\widehat{v_{r,k}^{n+1}|_{\Gamma(t^n)}}-\widehat{v_{r,k}^{n+1/2}|_{\Gamma(t^n)}}}{\Delta t}  \varphi_{r,k} d x
 \end{equation*}
 \begin{equation}
 = \int_0^R p_{in}(t^{n+1}) \varphi_{z,k}|_{z=0} dr - \int_0^R p_{out}(t^{n+1}) \varphi_{z,k}|_{z=L} dr, \label{step2d}
 \end{equation}
 where the initial data for the trace of the fluid velocity on $\Gamma(t^n)$ comes from the just calculated velocity of the
 thin structure $\widehat{v^{n+1/2}_{r,k}|_{\Gamma(t^n)}}=\hat{\xi}^{n+1/2}_{r,k} $.
  
  Then set $\hat{\boldsymbol U}_k^{n+1}=\hat{\boldsymbol U}_k^{n+1/2}$,  $\hat{\boldsymbol V}_k^{n+1}=\hat{\boldsymbol V}_k^{n+1/2}$,
$\hat{\eta}^{n+1}_{r,k}=\hat{\eta}^{n+1/2}_{r,k}$, and $\hat{\xi}^{n+1}_{r,k}=\hat{\xi}^{n+1/2}_{r,k}$.
 
 The ALE velocity $\boldsymbol{w}^{n+1}$ that appears in Problem {\bf A2} is calculated after computing the new  structure displacement in 
 Problem {\bf A1}. Note that Problem {\bf A1} is the only step in which the displacement of the structure changes. 
 Hence, after computing the new structure displacement in Problem {\bf A1}, $\hat{\eta}_{r,k}^n$ determines the ``reference domain'' and $\hat{\eta}_{r,k}^{n+1/2}$, which is the same as $\hat{\eta}_{r,k}^{n+1}$, determines the location of the new domain. 
To update the position of the fluid domain (after Problem {\bf A2}), we consider the following simple ALE mapping:
$$
A_{t^n}: \hat\Omega^f \to \Omega^f(t^n),\quad A_{t^n}(\hat{z},\hat{r}):=\left(\hat{z},\frac{R+\hat{\eta}_{r,k}^n}{R}\hat{r}\right)^\tau.
$$
We will  also need the explicit form of the ALE mapping from the computational 
domain $\Omega^f(t^n)$ to $\Omega^f(t^{n+1})$, which is given by
$$
A_{t^{n+1}} \circ A_{t^n}^{-1}: \Omega^f(t^n) \to \Omega^f(t^{n+1}), \quad A_{t^{n+1}} \circ A_{t^n}^{-1}(z,r) = \left(z, \frac{R + \hat{\eta}_{r,k}^{n+1}}{R+\hat{\eta}_{r,k}^n} r\right)^\tau.
$$

The corresponding Jacobian and the ALE velocity are given, respectively, by
\begin{equation}\label{jacobian}
J_n^{n+1} := \frac{R + \hat{\eta}_{r,k}^{n+1}}{R+\hat{\eta}_{r,k}^n}, \quad 
{\boldsymbol w}_k^{n+1}
=\frac{1}{\Delta t}\frac{\hat{\eta}_{r,k}^{n+1}-\hat{\eta}_{r,k}^n}{R+\hat{\eta}_{r,k}^n} \ {r}{\boldsymbol e}_r.
\end{equation}
Therefore, 
the time-derivative of the interface displacement is approximated by $(\hat{\eta}_{r,k}^{n+1}-\hat{\eta}_{r,k}^n)/\Delta t$, which enters 
the expression for the ALE velocity ${\boldsymbol w}^{n+1}$. 

\subsection{Stability analysis}
Let $\mathcal{E}_f^n$ denote the discrete energy of the fluid  problem, $\mathcal{E}_s^n$ the discrete energy of the structure problem, and let $\mathcal{E}_m^n$ denote the discrete energy of the simplified membrane problem at time level $n \Delta t$:
\begin{eqnarray*}
 \mathcal{E}_f^n & :=& \frac{\rho_f}{2} ||\boldsymbol v_k^n||^2_{L^2(\Omega^f(t^n))}, \\
 \mathcal{E}_s^n &:=& \frac{\rho_{s}}{2} ||\hat{\boldsymbol V}_k^n||^2_{L^2(\hat{\Omega}^s)} + \mu_s ||D(\hat{\boldsymbol U}_k^n)||^2_{L^2(\hat{\Omega^s})} + \frac{\lambda_s}{2}||\nabla \cdot \hat{\boldsymbol U}_k^n||^2_{L^2(\hat{\Omega}^s)}+ \frac{\gamma}{2} ||\hat{\boldsymbol U}_k^n||^2_{L^2(\hat{\Omega}^s)}, \\
 \mathcal{E}_m^n &:=&  \frac{\rho_{m}h}{2} ||\hat{\xi}_{r,k}^n||^2_{L^2(0,L)} + \frac{C_0}{2}||\hat{\eta}_{r,k}^n||^2_{L^2(0,L)}.
 \end{eqnarray*}
The following energy estimate holds for the full nonlinear FSI problem,
satisfying assumptions {\bf 1} and {\bf 2} above.

\begin{theorem}{\bf (An energy estimate of the operator splitting scheme)}
Let $\{(\boldsymbol v^n,\boldsymbol U^n, \boldsymbol V^n, \eta_r^n, \xi_r^n \}_{0 \le n \le N}$ be a solution of~\eqref{step1d}-\eqref{step2d}.
Then, the following energy estimate holds:
\begin{equation*}
\mathcal{E}_f^{N}+\mathcal{E}_s^{N}+\mathcal{E}_m^{N}+ \frac{\rho_f}{2}\sum_{n=0}^{N-1}
||\boldsymbol v_k^{n+1}-\boldsymbol v_k^{n+1/2}||^2_{L^2(\Omega^f(t^{n}))}
 + \mu_f \Delta t\sum_{n=0}^{N-1} ||D(\boldsymbol v_k^{n+1})||^2_{L^2(\Omega^f(t^{n}))} 
 \end{equation*}
\begin{equation*}
 + \frac{\rho_m h}{2}\sum_{n=0}^{N-1}
||\widehat{v_{r,k}^{n+1}|_{\Gamma(t^n)}}-\widehat{v_{r,k}^{n+1/2}|_{\Gamma(t^n)}}||^2_{L^2(0,L)}
\end{equation*}
\begin{equation*}
\le \mathcal{E}_f^0 + \mathcal{E}_s^0+ \mathcal{E}_m^0 +  \frac{C}{\mu_f}\Delta t \sum_{n=0}^{N-1}  ||p_{in}(t^n)||^2_{L^2(0,R)}+  \frac{C}{\mu_f}\Delta t \sum_{n=0}^{N-1}  ||p_{out}(t^n)||^2_{L^2(0,R)}.
\end{equation*}
 \end{theorem}
 
 This energy estimate shows that the total energy of the coupled discretized problem is bounded by the energy of the initial data,
 and by the work done by the inlet and outlet dynamic pressure data, independently of the time step, for all the parameters
 in the problem. Notice that the total energy of the problem includes
 the kinetic energy of the fluid and the kinetic energy of both structures, the elastic energy of both structures, and, additionally,
 the kinetic energy that is due to the motion of the fluid domain and of the fluid domain boundary. The last two contributions are accounted for
 in the two terms on the left hand-side of the 
 energy inequality, containing the $L^2$-norms of the difference between the fluid velocities at times $n+1$ and $n+1/2$. 

 \begin{proof}
To prove the energy estimate, we test the structure  problem~\eqref{step1d} with 
$
(\hat{\boldsymbol \psi}_k, \hat{\boldsymbol \phi}_k)=(\displaystyle\frac{\hat{\boldsymbol U}_k^{n+1/2}-\hat{\boldsymbol U}_k^{n}}{\Delta t},\displaystyle\frac{\hat{\boldsymbol V}_k^{n+1/2}-\hat{\boldsymbol V}_k^{n}}{\Delta t}),
$
and the fluid problem~\eqref{step2d} with
$
(\boldsymbol \varphi_k, q_k) = (\boldsymbol v_k^{n+1}, p_k^{n+1}).
$
With this choice of test functions, after multiplying by $\Delta t$, Problem {\bf A1} 
 reduces to
\begin{equation*}
\mathcal{E}_s^{n+1/2}+\frac{\rho_{m}h}{2} ||\hat{\xi}_{r,k}^{n+1/2}||^2_{L^2(0,L)}-\frac{\rho_{m}h}{2} ||\hat{\xi}_{r,k}^n||^2_{L^2(0,L)}+ \frac{C_0}{2}||\hat{\eta}_{r,k}^{n+1/2}||^2_{L^2(0,L)}
\end{equation*}
\begin{equation}
 - \frac{C_0}{2}||\hat{\eta}_{r,k}^n||^2_{L^2(0,L)} = \mathcal{E}_s^{n}. \label{en1} 
\end{equation}

In Problem {\bf A2}, after replacing the test functions by the fluid velocity $\boldsymbol v_k^{n+1}$, 
the symmetrized advection terms cancel out.
To deal with the first term on the left in Problem {\bf A2}, we use the identity
\begin{equation}\label{identity}
a(a-b) =\frac 1 2( a^2-b^2+(a-b)^2).
\end{equation}
Then, the first two terms  combined (fluid inertia) become
$$
\frac {\rho_f}{2}\frac{1}{\Delta t}\int_{\Omega^f(t^n)}\left(\bigg(1+\frac{\hat{\eta}_{r,k}^{n+1}-\hat{\eta}_{r,k}^n}{R+\hat{\eta}_{r,k}^n}\bigg) |{\boldsymbol v}_k^{n+1}|^2
+|{\boldsymbol v}_k^{n+1}-{\boldsymbol v}_k^{n+1/2}|^2-|{\boldsymbol v}_k^{n+1/2}|^2 \right) d \boldsymbol x
$$
$$
=\frac{\rho_f}{2}\frac{1}{\Delta t}\int_{\Omega^f(t^n)}\Big (\frac{R+\hat{\eta}_{r,k}^{n+1}}{R+\hat{\eta}_{r,k}^n}|{\boldsymbol v}_k^{n+1}|^2
+|{\boldsymbol v}_k^{n+1}-{\boldsymbol v}_k^{n+1/2}|^2-|{\boldsymbol v}_k^{n+1/2}|^2\Big)d \boldsymbol x .
$$
Note that, since the fluid does not change in Problem {\bf A1}, $\boldsymbol v^{n+1/2}_k = \boldsymbol v_k^n$. Furthermore,
notice that $({R+\eta^{n+1}})/({R+\eta^n})$ is exactly the Jacobian of the ALE mapping from $\Omega^f(t^n)$ to $\Omega^f(t^{n+1})$, 
see \eqref{jacobian}, and so we can convert that integral into an integral over $\Omega^f(t^{n+1})$ to recover
the kinetic energy of the fluid at the next time-step:
$$
\frac{\rho_f}{2}\frac{1}{\Delta t}\int_{\Omega^f(t^n)}\frac{R+\hat{\eta}_{r,k}^{n+1}}{R+\hat{\eta}_{r,k}^n}|{\boldsymbol v}_k^{n+1}|^2 d \boldsymbol x
=\frac{\rho_f}{2}\frac{1}{\Delta t}\int_{\Omega^f(t^{n+1})}|{\boldsymbol v}_k^{n+1}|^2 d \boldsymbol x.
$$
This calculation also shows that the ALE mapping and its Jacobian satisfy the geometric conservation law property, as studied 
by Farhat et al. in \cite{Farhat}. 

Therefore, after multiplying by $\Delta t$, the energy associated with Problem {\bf A2} is given by 
\begin{equation*}
\frac{\rho_f}{2}||\boldsymbol v_k^{n+1}||^2_{L^2(\Omega^f(t^{n+1}))}-\frac{\rho_f}{2}||\boldsymbol v_k^{n+1/2}||^2_{L^2(\Omega^f(t^{n}))}+\frac{\rho_f}{2}||\boldsymbol v_k^{n+1}-\boldsymbol v_k^{n+1/2}||^2_{L^2(\Omega^f(t^{n}))} 
\end{equation*}
\begin{equation*}
+2 \mu_f \Delta t ||D(\boldsymbol v_k^{n+1})||^2_{L^2(\Omega^f(t^{n}))}+\frac{\rho_m h}{2}||\widehat{ v_{r,k}^{n+1}|_{\Gamma(t^n)}}||^2_{L^2(0,L)} -\frac{\rho_m h}{2}||\widehat{v_{r,k}^{n+1/2}|_{\Gamma(t^n)}}||^2_{L^2(0,L))}
\end{equation*}
\begin{equation*}
+\frac{\rho_m h}{2}||\widehat{v_{r,k}^{n+1}|_{\Gamma(t^n)}}-\widehat{v_{r,k}^{n+1/2}|_{\Gamma(t^n)}}||^2_{L^2(0,L)}=\Delta t  \int_0^R p_{in}(t^{n+1}) v_{z,k}^{n+1}|_{z=0} dr
\end{equation*}
\begin{equation} 
 - \Delta t \int_0^R p_{out}(t^{n+1}) v_{z,k}^{n+1}|_{z=L} dr. \label{en3}
\end{equation}

Now, equations ~\eqref{en1} and~\eqref{en3} are added, and the following properties are taken into account:
first $\widehat{v_{r,k}^{n+1/2}|_{\Gamma(t^n)}}=\hat{\xi}_{r,k}^{n+1/2}$, and further,
due to the kinematic coupling condition 
we also have $\widehat{ v_{r,k}^{n+1}|_{\Gamma(t^n)}} = \hat{\xi}^{n+1}_{r,k}$.
After adding equations~\eqref{en1} and~\eqref{en3}, and taking into account $\mathcal{E}_s^{n+1/2} = \mathcal{E}_s^{n+1}$  since the displacement changes only in the first step, we get
\begin{equation*}
\mathcal{E}_f^{n+1}+\mathcal{E}_s^{n+1}+\mathcal{E}_m^{n+1}+ \frac{\rho_f}{2}
||\boldsymbol v_k^{n+1}-\boldsymbol v_k^{n+1/2}||^2_{L^2(\Omega^f(t^{n}))}
 +2 \mu_f \Delta t ||D(\boldsymbol v_k^{n+1})||^2_{L^2(\Omega^f(t^{n}))} 
 \end{equation*}
\begin{equation*}
 + \frac{\rho_m h}{2}
||\widehat{v_{r,k}^{n+1}|_{\Gamma(t^n)}}-\widehat{v_{r,k}^{n+1/2}|_{\Gamma(t^n)}}||^2_{L^2(0,L)}
\end{equation*}
\begin{equation*}
 = \mathcal{E}_f^n + \mathcal{E}_s^n+ \mathcal{E}_m^n +  \Delta t  \int_0^R p_{in}(t^{n+1}) v_{z,k}^{n+1}|_{z=0} dr 
 - \Delta t \int_0^R p_{out}(t^{n+1}) v_{z,k}^{n+1}|_{z=L} dr.
\end{equation*}

To bound the right-hand side of this equality, we use the Cauchy-Schwartz and Young's inequalities
\begin{equation*}
\Delta t  \int_0^R p_{in}(t^{n+1}) v_{z,k}^{n+1}|_{z=0} dr  - \Delta t \int_0^R p_{out}(t^{n+1})  v_{z,k}^{n+1}|_{z=L} dr 
\end{equation*}
\begin{equation*}
\le \frac{\Delta t}{2 \epsilon_1 } ||p_{in}(t^n)||^2_{L^2(0,R)}+\frac{\Delta t}{2 \epsilon_1 } ||p_{out}(t^n)||^2_{L^2(0,R)}+\epsilon_1 \Delta t
|| {\boldsymbol v}_k^{n+1}||^2_{L^2(0,R)}.
\end{equation*}
By the trace and Korn inequalities, we then have
\begin{equation*}
|| {\boldsymbol v}_k^{n+1}||^2_{L^2(0,R)} \le  C ||\boldsymbol D( \boldsymbol v_k^{n+1})||^2_{L^2(\Omega^f(t^n))},
\end{equation*}
where $C$ is the constant from the trace and Korn inequalities. 
In general, Korn's constant depends on the domain. It was shown, however, that for domains associated 
with fluid-structure interaction problems of the type studied in this manuscript, the Korn's constant is independent
of the sequence of approximating domains \cite{BorSun,BorSunMulti}.

By setting $\epsilon_1 = \displaystyle\frac{\mu_f}{C}$, the last term can be combined with the term 
on the left hand-side, associated with the fluid diffusion. Finally, summing the inequality from $n=0$ to $N-1$ we prove the desired energy estimate. 
\end{proof}

Notice that our operator splitting algorithm consists of solving a set of sub-problems which are linear (elliptic). 
It was shown in \cite{BorSunMulti} using compactness arguments, that the sequence of linear sub-problems converges to a weak solution 
of the full,
nonlinear FSI problem with multilayered structures, in which the thin structure was modeled by the linear wave equation. 
The energy estimate, presented above, shows that the combined linear discretized sub-problems are stable, and that
the sequence of approximate solutions to the nonlinear FSI problem is uniformly bounded in the corresponding energy norms.

Numerical results are presented next, where two instructive numerical examples are solved using the proposed scheme.
The first example presents a simplified FSI problem with two layers, for which an exact solution was found and compared
with the numerically calculated solution. The second example concerns a full, nonlinear FSI problem with two structural layers
in which the thickness of the thin layer is taken to converge to zero. It was shown that the final solution compares well with the solution of
the limiting FSI problem, obtained using a different computational solver from the one presented in this manuscript.
Both examples show that the scheme is unconditionally stable, and that it converges to a solution of the underlying problem. 



\section{Numerical example 1: A simplified FSI problem with exact solution.}\label{sec:numerics1}

We consider a simplified FSI problem with multiple structural layers that satisfies the following simplifying assumptions:
\begin{itemize}
 \item[1.] The fluid problem is defined on the fixed, reference domain of width $R$, and length $L$ (the coupling is linear).
 \item[2.] The fluid problem is driven by the constant inlet and outlet pressure data $p_{in}$ and $p_{out}= 0$ (the pressure drop is constant).
 \item[3.] Only radial displacement of the thin and thick structure is assumed to be different from zero.
\end{itemize}
Assumption 3 implies that
the thin structure membrane model takes the form:
\begin{equation*}
 \rho_K h \frac{\partial^2 \eta_r}{\partial t} + C_0 \eta_r = f_r,
\end{equation*}
while the thick structure problem simplifies as follows:
\begin{equation*}
 \rho_s \frac{\partial^2 d_r}{\partial t^2} = \mu \frac{\partial^2 d_r}{\partial x^2} + (\mu+\lambda)\frac{\partial^2 d_r}{\partial y^2}.
\end{equation*}
Finally, the coupling conditions between the fluid and the multilayered structure are given by
\begin{align*}
& f_r = p + (\lambda+\mu) \frac{\partial d_r}{\partial y} & \textrm{on} \; \Gamma \times(0,T), \\
&\frac{\partial \eta_r}{\partial t} = u_r & \textrm{on} \; \Gamma \times(0,T), \\
&\eta_r = U_r & \textrm{on} \; \Gamma \times(0,T).
 \end{align*}
The exact solution to this problem is given by the following. The fluid flow through the fixed cylinder with constant pressure drop 
is given by the Poiseuille velocity profile: 
\begin{equation*}
 u^{e}_z (z,r) = u^e_z(r)= \frac{p_{in}-p_{out}}{2 \mu_F L}(R^2-r^2), \quad u^{e}_r=0, 
\end{equation*}
and the fluid pressure is linear within the channel:
\begin{equation*}
p^e(z,r)=p^e(z) = \frac{p_{out}z+p_{in}(L-z)}{L}, \ z\in(0,L), \ r\in(0,R). 
\end{equation*}
The radial displacements of the thin and thick structure are given by:
\begin{equation*}
 \eta_r^e (z)= \frac{p^e(z)}{C_0}, \quad d^e_r(z,r) = d_r^e(z) = \eta_r^e(z).  
\end{equation*}

{{
\begin{center}
\begin{table}[ht!]
{\small{
\begin{tabular}{|l l l l |}
\hline
\textbf{Parameters} & \textbf{Values} & \textbf{Parameters} & \textbf{Values}  \\
\hline
\hline
\textbf{Radius} $R$ (cm)  & $0.5$  & \textbf{Length} $L$ (cm) & $6$  \\
\hline
\textbf{In. press.} $p_{in}$ (dyne/cm$^2$)& $250$ &\textbf{Out. press.} $p_{out}$ (dyne/cm$^2$) & $0$    \\
\textbf{Fluid density} $\rho_f$ (g/cm$^3$)& $1$ &\textbf{Dyn. viscosity} $\mu$ (g/cm s) & $0.35$    \\
\hline
\textbf{Thin wall:}  &  &  &  \\
\textbf{Density} $\rho_{m} $(g/cm$^3$) & $1.1$  & \textbf{Thickness} $h$ (cm) & $0.02$  \\
\textbf{Lam\'e coeff.} $\mu_{m} $(dyne/cm$^2$) & $1.07 \times 10^6$  & \textbf{Lam\'e coeff.} $\lambda_{m} $(dyne/cm$^2$) & $4.29 \times 10^6$ \\
\hline
\textbf{Thick wall:}  &  &  &  \\
\textbf{Density} $\rho_{s} $(g/cm$^3$) & $1.1$  & \textbf{Thickness} $H$ (cm) & $0.1$  \\
\textbf{Lam\'e coeff.} $\mu_{s} $(dyne/cm$^2$) & $1.07 \times 10^6$  & \textbf{Lam\'e coeff.} $\lambda_{s} $(dyne/cm$^2$) & $4.29 \times 10^6$ \\
\textbf{Spring coeff.} $\gamma $(dyne/cm$^4$) & $0$  &  &    \\
\hline
\end{tabular}
}}\caption{Geometry, fluid and structure parameters used in Example 1.}
\label{T0}
\end{table}
\end{center}
}}
We solve this problem numerically using the parameters given in Table~\ref{T0}.
The initial data was taken to be 
$$
{\bf u} = 0, p = p_{out}, \eta_r = 0, d_r = 0, \ {\rm at}\  t = 0,
$$
while at the inlet and outlet boundaries we kept both structures fixed, with the inlet and outlet displacement data 
tailored so that the final solution does not exhibit a boundary layer:
$$\eta_r|_{z=0} = d_r|_{z=0} = \frac{p_{in}}{C_0}, \quad \eta_r|_{z=L} = d_r|_{z=L} = \frac{p_{out}}{C_0}  = 0,\forall t > 0.$$

The numerical scheme with $\beta = 1$ was implemented, and the problem was solved until 
the steady state was achieved. 
With the time step  $\Delta t = 10^{-5}$ it took 200 iterations to achieve  the accuracy of less than 0.08\%.
Namely, the maximum relative error between the computed and exact solution was less than 0.08\% (namely, 0.000778).

 \begin{figure}[ht!]
 \centering{
 \includegraphics[scale=0.70]{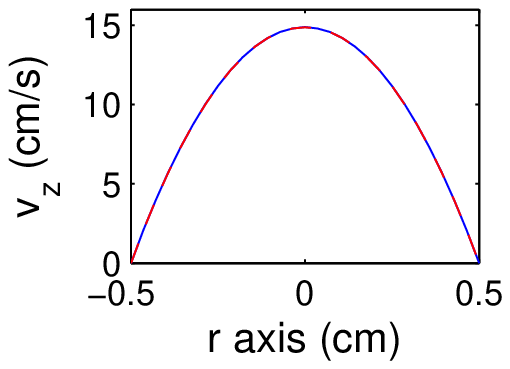}
  \includegraphics[scale=0.70]{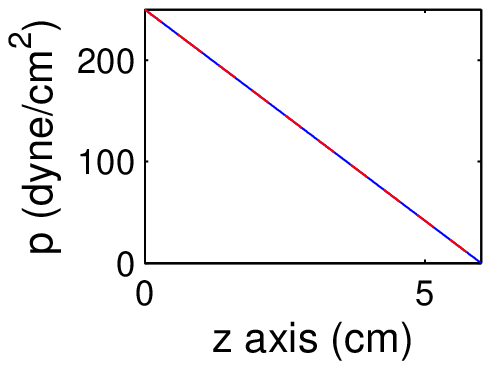}
   \includegraphics[scale=0.70]{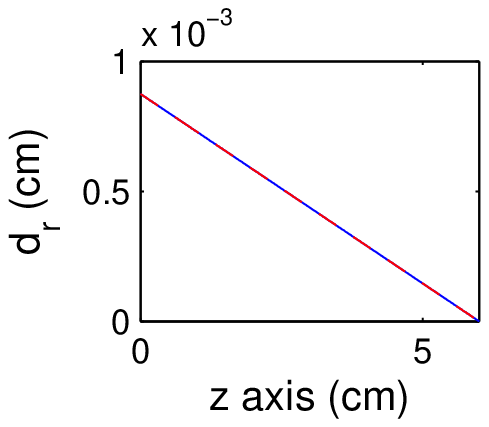}
 }
 \caption{Comparison between the computed solution (in blue) and the exact solution (in red). 
 The two are superimposed. Left: Axial velocity. Middle: Fluid pressure. Right: Radial displacement.}
\label{Ex1F1}
 \end{figure}

Figure~\ref{Ex1F1} shows a comparison between the computed (blue) and the exact solution (red) for axial velocity (left), fluid pressure (middle), 
and radial displacement (right), showing excellent agreement.
The corresponding relative errors are given by the following:
$$\frac{||{\bf u}^e - {\bf u}||_{L^2(\Omega^f)}}{||{\bf u}^e||_{L^2(\Omega^f)}} = 7.78 \times 10^{-4}, \quad \frac{||p^e - p||_{L^2(\Omega^f)}}{||p^e||_{L^2(\Omega^f)}} = 1.17 \times 10^{-4},$$
$$\frac{|| \eta_r^e - \eta_r||_{L^2(0,L)}}{||\eta_r^e||_{L^2(0,L)}} = 3.82 \times 10^{-5}, \quad\frac{|| d_r^e - d_r||_{L^2(\Omega^s)}}{||d_r^e||_{L^2(\Omega^s)}} = 3.82 \times 10^{-5}.$$
We conclude that the scheme behaves well for this simplified FSI problem with multiple structural layers. 

\section{Numerical example 2: A fully nonlinear FSI problem.}\label{sec:numerics2}
In this example we solve the full, nonlinear FSI problem \eqref{NS1}-\eqref{dynamic}.
In this example the structure consists of two layers, one with thickness $h$ (thin) and one with thickness $H>h$ (thick).
The combined thickness of the composite structure is $h + H = 0.12$cm. 
The data and parameters used in the simulation correspond to those used in  \cite{thick} and \cite{badia}
to test different FSI solvers. They are given in Table~\ref{T1}.
{{
\begin{center}
\begin{table}[ht!]
{\small{
\begin{tabular}{|l l l l |}
\hline
\textbf{Parameters} & \textbf{Values} & \textbf{Parameters} & \textbf{Values}  \\
\hline
\hline
\textbf{Radius} $R$ (cm)  & $0.5$  & \textbf{Length} $L$ (cm) & $6$  \\
\hline
\textbf{Fluid density} $\rho_f$ (g/cm$^3$)& $1$ &\textbf{Dyn. viscosity} $\mu$ (g/cm s) & $0.035$    \\
\hline
\textbf{Thin wall:}  &  &  &  \\
\textbf{Density} $\rho_{m} $(g/cm$^3$) & $1.1$  & \textbf{Thickness} $h$ (cm) & $0.02$  \\
\textbf{Lam\'e coeff.} $\mu_{m} $(dyne/cm$^2$) & $5.75 \times 10^5$  & \textbf{Lam\'e coeff.} $\lambda_{m} $(dyne/cm$^2$) & $1.7 \times 10^6$ \\
\hline
\textbf{Thick wall:}  &  &  &  \\
\textbf{Density} $\rho_{s} $(g/cm$^3$) & $1.1$  & \textbf{Thickness} $H$ (cm) & $0.1$  \\
\textbf{Lam\'e coeff.} $\mu_{s} $(dyne/cm$^2$) & $5.75 \times 10^5$  & \textbf{Lam\'e coeff.} $\lambda_{s} $(dyne/cm$^2$) & $1.7 \times 10^6$ \\
\textbf{Spring coeff.} $\gamma $(dyne/cm$^4$) & $4 \times 10^6$  &  &    \\
\hline
\end{tabular}
}}\caption{Geometry, fluid and structure parameters that are used in Example 2.}
\label{T1}
\end{table}
\end{center}
}}
Since there are no benchmark results in FSI literature containing a two-layered composite structure interacting with an incompressible, viscous fluid,
we take the following approach to test our FSI solver. We solve 
a sequence of FSI problems \eqref{NS1}-\eqref{dynamic} in which the thickness
of the thin layer converges to zero, $h \to 0$, while the combined thickness $h + H$ remains fixed, and equal to $0.12$cm. 
First, we proved, analytically, that the limiting functions (displacement and velocity) satisfy the FSI problem containing a single, thick structure,
with thickness $0.12$cm. This is presented in the Appendix. 
Then we numerically compare the ``limiting'' solution of the sequence of  FSI problems with composite structures, 
obtained using the FSI solver proposed in this manuscript, with a solution 
of the FSI problem containing a single, thick structure with thickness $h+H = 0.12$cm, obtained using a different computational solver,
namely a solver developed in \cite{thick}. We show that the two are in excellent agreement.

The elastodynamics of the thin structural layer is modeled using the linearly elastic Koiter
membrane equations with both radial and longitudinal displacement \eqref{structure1}, \eqref{structure2},
while the elastodynamics of the thick structure is modeled using the equations of 2D linear elasticity \eqref{linela}.

The flow is driven by the time-dependent pressure data:
\begin{equation*}\label{pressure}
 p_{in}(t) = \left\{\begin{array}{l@{\ } l} 
\frac{p_{max}}{2} \big[ 1-\cos\big( \frac{2 \pi t}{t_{max}}\big)\big] & \textrm{if} \; t \le t_{max}\\
0 & \textrm{if} \; t>t_{max}
 \end{array} \right.,   \quad p_{out}(t) = 0 \;\forall t \in (0, T),
\end{equation*}
where $p_{max} = 1.333 \times 10^4$ (dyne/cm$^2$) and $t_{max} = 0.003$ (s). 

 \begin{figure}[ht!]
 \centering{
 \includegraphics[scale=0.56]{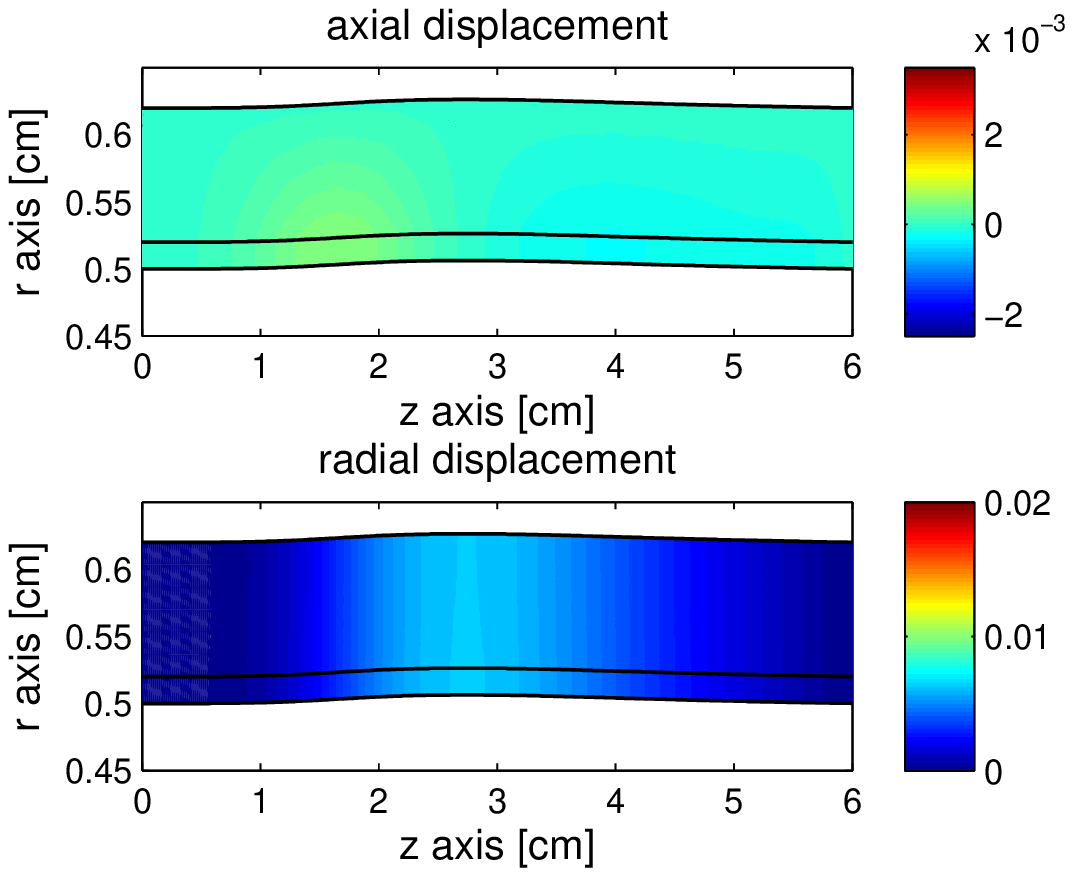}
 \includegraphics[scale=0.56]{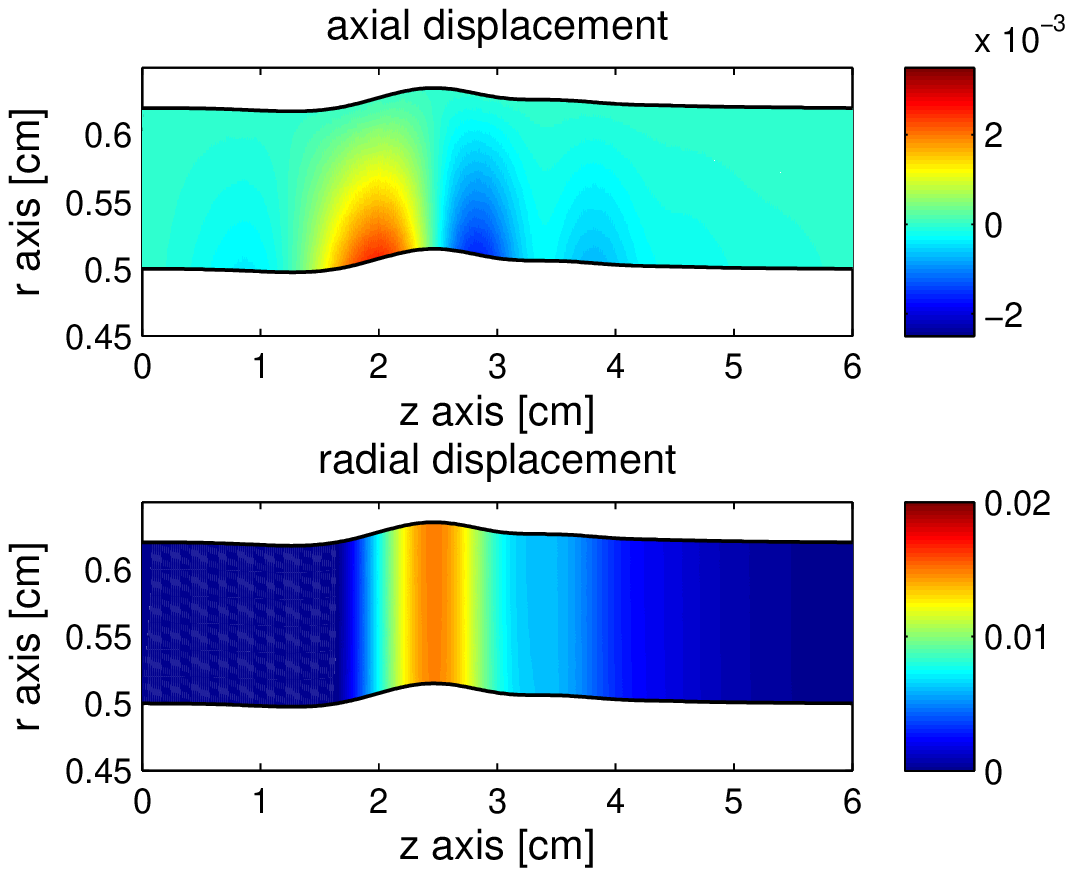}
 }
 \caption{Axial displacement (top) and radial displacement (bottom) at time $t = 8$ ms obtained using the model 
 capturing two structural layers (left), and the model capturing FSI with a single thick structural layer \cite{thick}  (right).}
\label{Ex2F1}
 \end{figure}
 
  \begin{figure}[ht!]
 \centering{
 \includegraphics[scale=0.35]{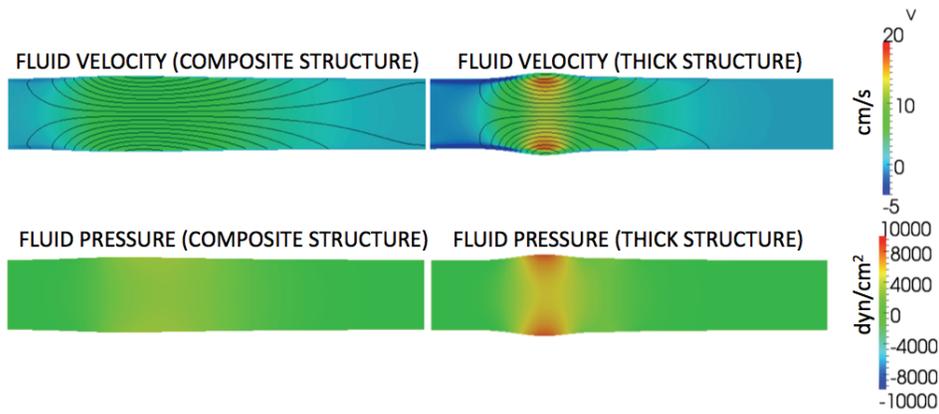}
 }
 \caption{Fluid velocity (top) and fluid pressure (bottom) at time $t = 8$ ms obtained using the model 
 capturing two structural layers (left), and the model capturing FSI with a single thick structural layer \cite{thick}  (right).}
\label{PandVel}
 \end{figure}
 
  \begin{figure}[ht!]
 \centering{
 \includegraphics[scale=0.7]{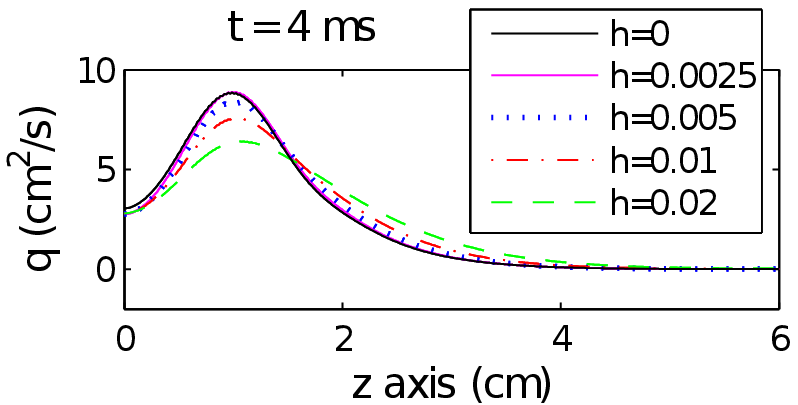}
  \includegraphics[scale=0.7]{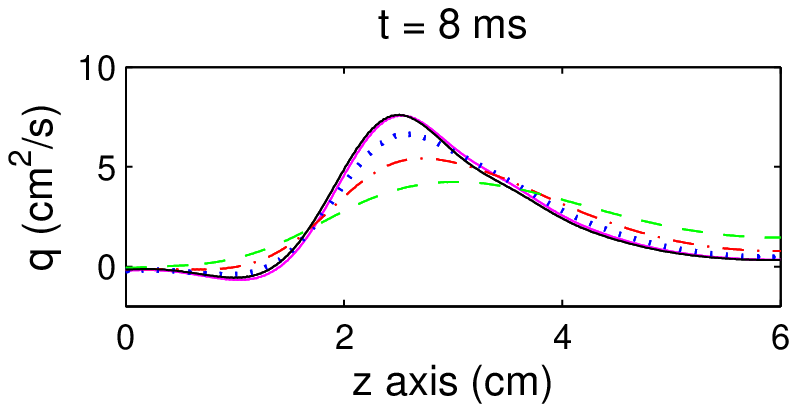}
   \includegraphics[scale=0.7]{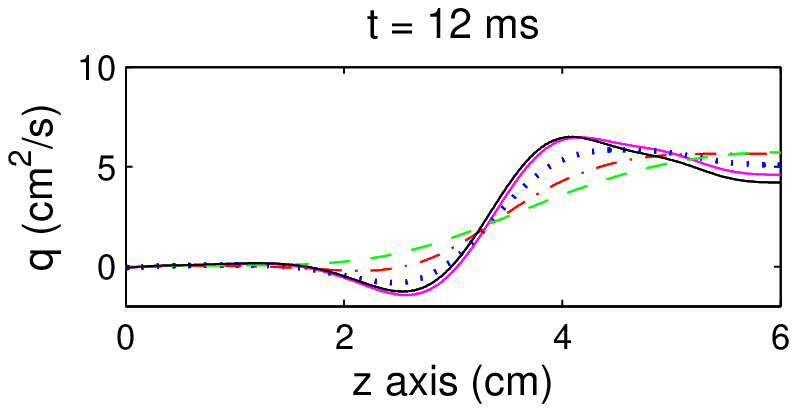}
 }
 \caption{Flowrate computed using two different models: the model in  \cite{thick} containing a single thick structural layer ($h=0$), and the model considered in this manuscript, consisting of two layers. The thickness of the thin membrane layer was decreased from $h = 0.02$ to $h = 0.0025$ cm.
 The combined thickness of the two-layered composite structure was kept constant at $h+H = 0.12$cm.
Convergence of the solutions to the 
FSI solution containing a single, thick structure ($h = 0$) can be observed.
}
\label{Ex2F3}
 \end{figure}
 
\begin{figure}[ht!]
 \centering{
 \includegraphics[scale=0.7]{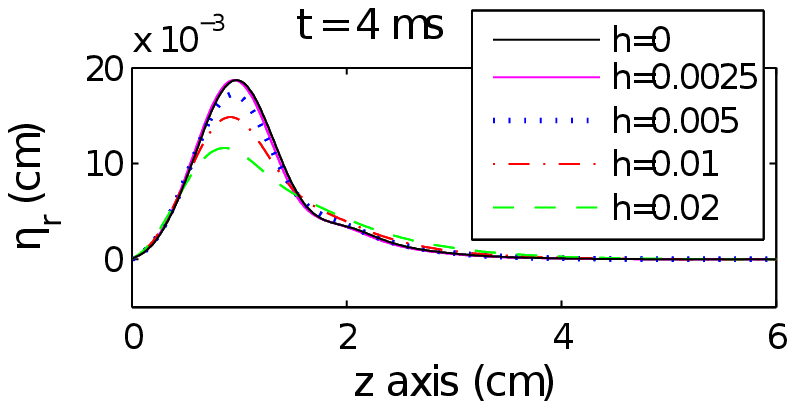}
  \includegraphics[scale=0.7]{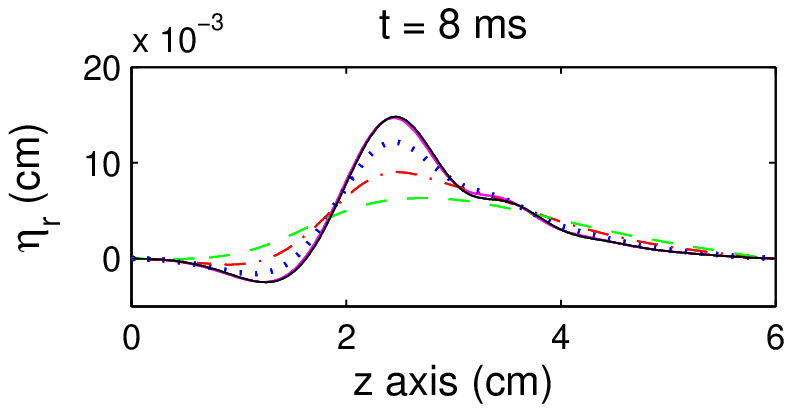}
   \includegraphics[scale=0.7]{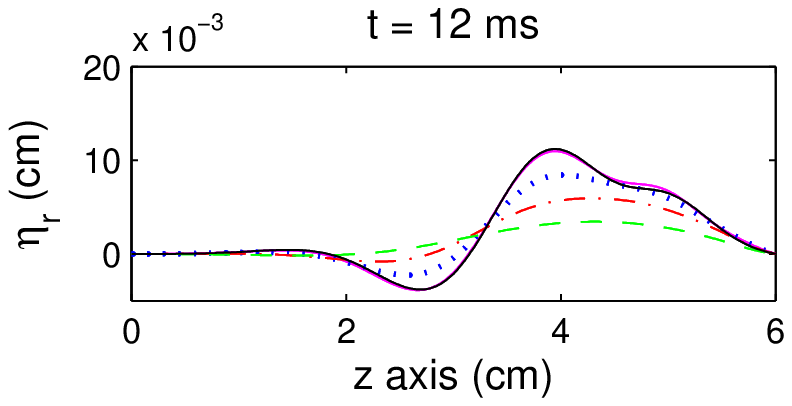}
 }
 \caption{Displacement of the fluid-structure interface obtained under the same conditions as those described in Figure~\ref{Ex2F3}. }
\label{Ex2F4}
 \end{figure}
The kinematically-coupled $\beta$ scheme, described in Sections~\ref{splitting} and \ref{discretization}, was used to solve the multi-layered FSI problem with $\beta = 1$,
while the scheme presented in \cite{thick} was used to solve the single-layered FSI benchmark problem with thick structure.
The problem was solved over the time interval $[0,0.012]$s, using the time step $\Delta t = 5 \times 10^{-5}$.
 Figure~\ref{Ex2F1} shows the axial and radial displacement at time $t = 8$ ms obtained using the multilayered model (left) and the single-layered model (right) for the arterial wall. Figure~\ref{PandVel} shows the corresponding fluid velocity and pressure.
 One can notice significant smoothing of both the displacement as well as the fluid velocity and pressure in the composite, i.e., multi-layered 
 structure case. Same data are used for both simulations.
  \begin{figure}[ht!]
 \centering{
 \includegraphics[scale=0.7]{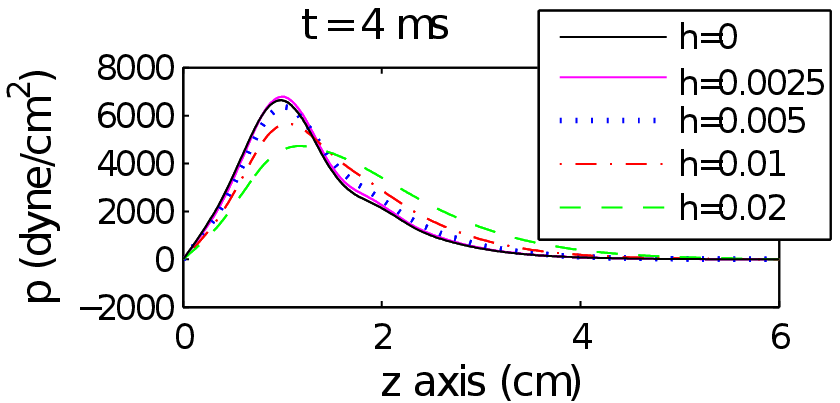}
  \includegraphics[scale=0.7]{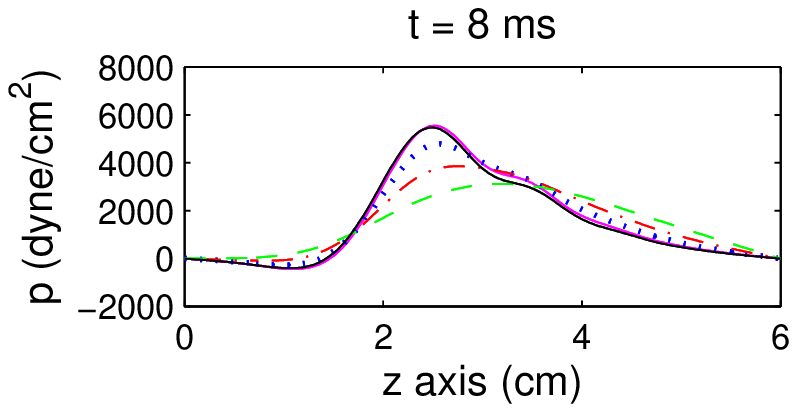}
   \includegraphics[scale=0.7]{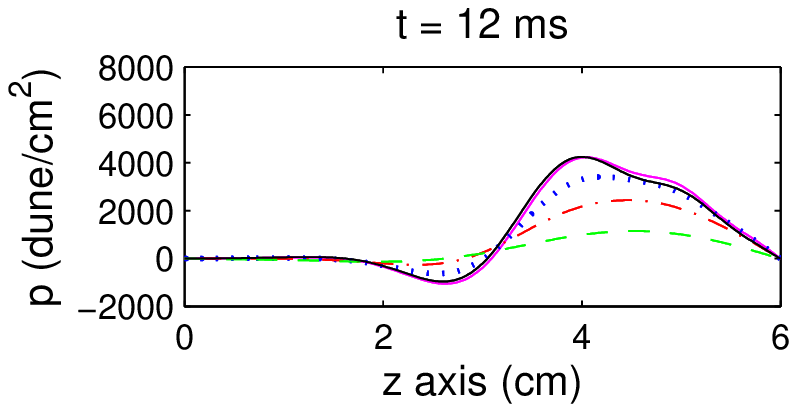}
 }
 \caption{Mean pressure obtained under the same conditions as those described in Figure~\ref{Ex2F3}.}
\label{Ex2F5}
 \end{figure}
 We further compared the results of the multilayered model with the single layered model as the thickness of the thin structure $h$ goes to zero. As we decreased $h$, we increased $H$ to maintain the constant combined thickness $h+H = 0.12$ cm. 
 Figures~\ref{Ex2F3},~\ref{Ex2F4} and~\ref{Ex2F5} show the flowrate, mean pressure and displacement of the fluid-structure interface obtained using different values of $h$. The results obtained using the single layered wall model correspond to the label $h=0$. Indeed, we can see that as we decrease the thickness of the fluid-structure interface, the numerical results obtained using our multilayered model, approach the results obtained using the single-layered FSI model! Notice how for $h=0.025$cm the solutions obtained using the multilayered model
 and the single thick structure model ($h = 0$ in Figures~\ref{Ex2F3}, \ref{Ex2F4}, \ref{Ex2F5}) are almost identical.

\begin{figure}[ht!]
 \centering{
 \includegraphics[scale=0.4]{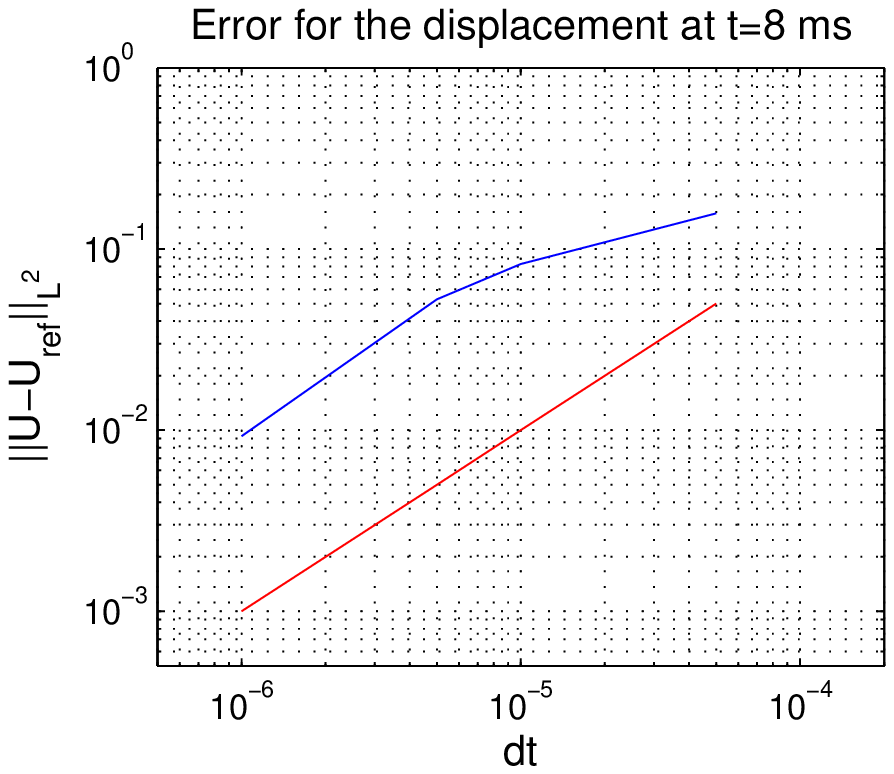}
  \includegraphics[scale=0.4]{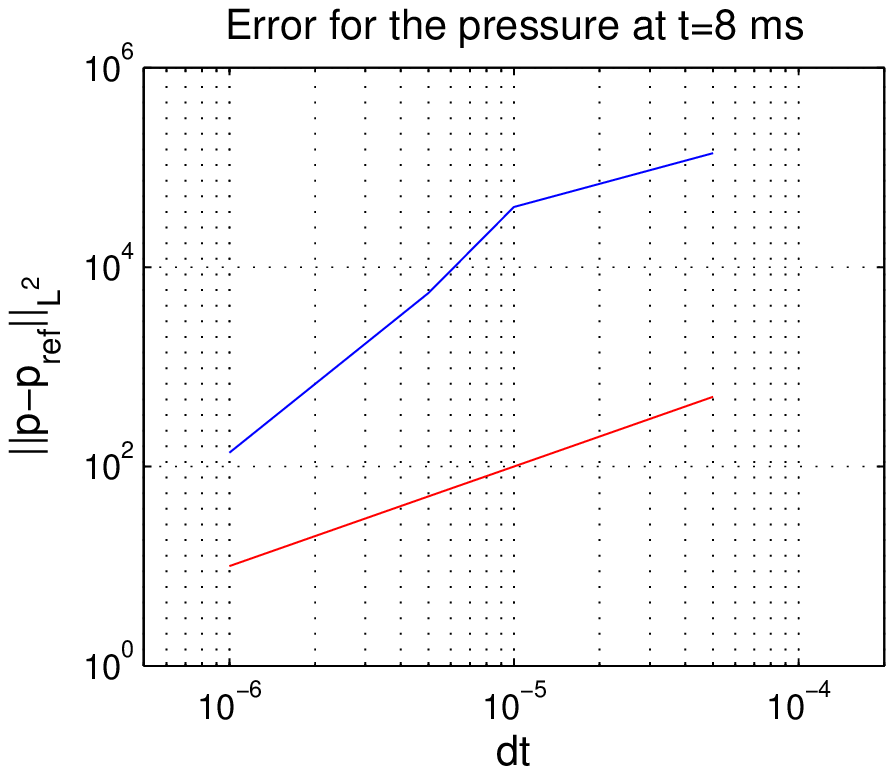}
   \includegraphics[scale=0.4]{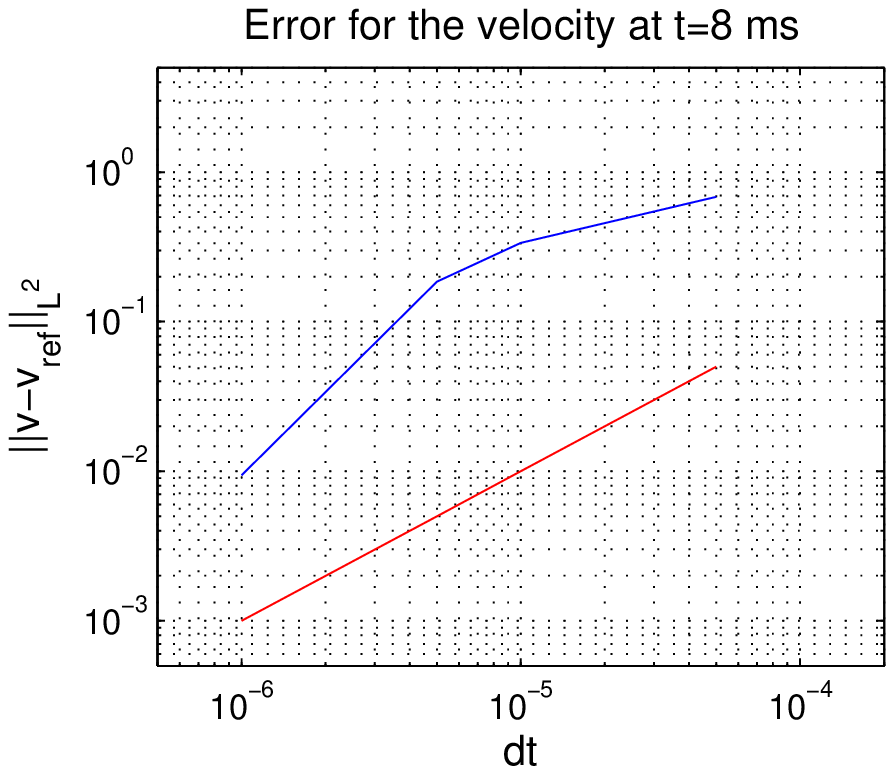}
 }
 \caption{The $L^2$-norm error of the structure displacement, fluid pressure and velocity vs. time step, shown in the log-log plots.
 The red line shows the slope associated with first-order accuracy. We see 1st-order accuracy for displacement, and higher-order 
 accuracy for pressure and velocity.}
\label{ERROR}
 \end{figure}

We conclude this section by showing the log-log plots of computational error versus time discretization step, presented in Figure~\ref{ERROR}. 
The red line in the plots indicates the slope associated with 1st-order accuracy. 
We see that our computational solver is 1st-order accurate in time for structure displacement calculations, while producing higher order temporal 
accuracy for fluid pressure and velocity.

\subsection{Regularizing Effects by Thin Fluid-Structure Interface with Mass}\label{reg}

We conclude this section with a remark on the regularizing effects of the thin fluid-structure interface with mass.
Figures~\ref{Ex2F3}, \ref{Ex2F4}, \ref{Ex2F5} indicate that as we increase 
inertia of the thin fluid-structure interface with mass by increasing its thickness, the solution of the entire FSI problem is 
damped, or regularized. More precisely, if one looks at the FSI problem with a single thick structural layer, the fluid-structure 
interface is simply the massless trace of the thick structure that is in contact with the fluid. Mathematically, in that case the trace
of the structure displacement is not well-defined (assuming regularity of 
the data consistent with weak solutions),
and using energy estimates it is not possible to even show that the fluid-structure interface is continuous. 
In the case when the fluid-structure interface has mass, it was shown in \cite{BorSunMulti} that not only is the fluid-structure interface
continuous, but its evolution can be controlled by the energy norm of the time derivative of its displacement. 
We see effects of this in the solutions presented in Figures~\ref{Ex2F3}, \ref{Ex2F4}, \ref{Ex2F5}, and in Figures~\ref{Ex2F2} below.
In Figure~\ref{Ex2F2} below we focus on the  displacement and  displacement velocity of the fluid-structure interface, which measures
the effects of inertia.
  \begin{figure}[ht!]
 \centering{
 \includegraphics[scale=0.48]{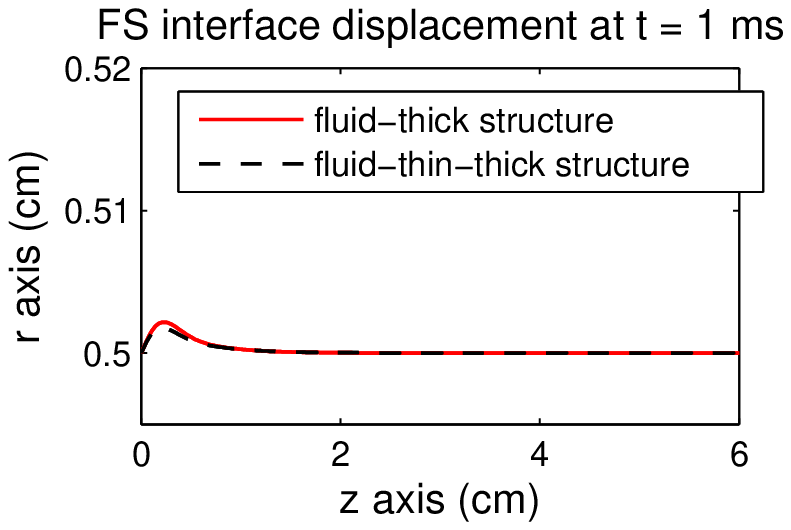}
   \includegraphics[scale=0.48]{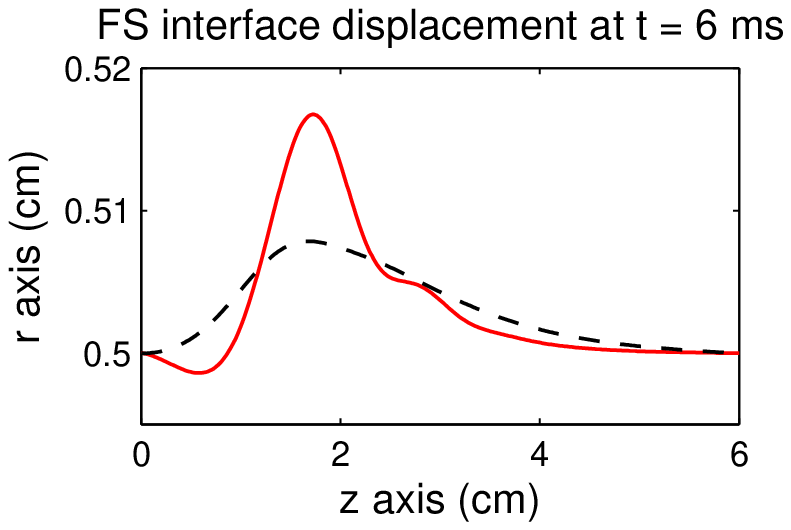}
     \includegraphics[scale=0.48]{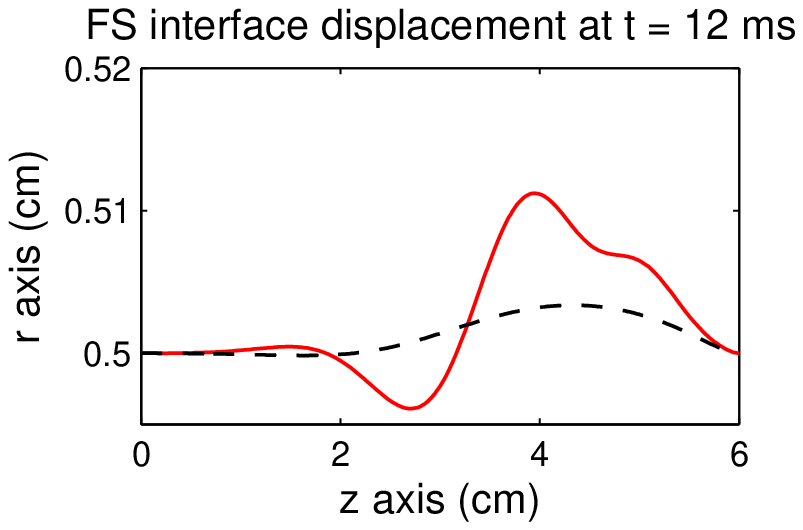}
  \includegraphics[scale=0.48]{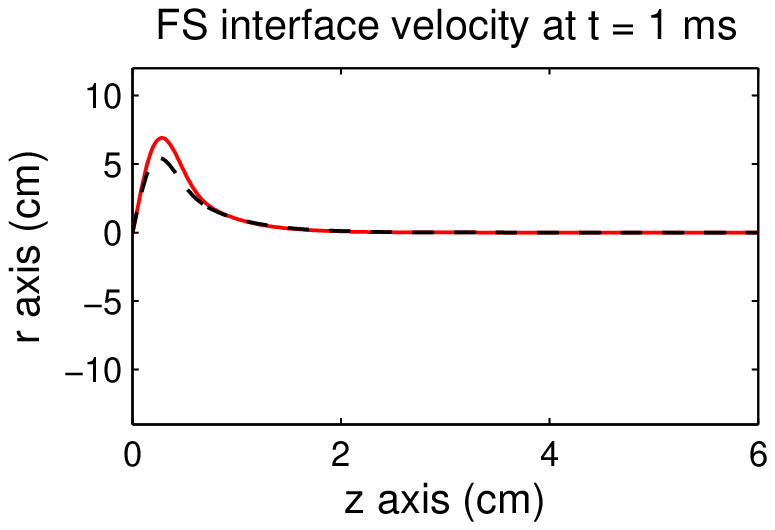}
  \includegraphics[scale=0.48]{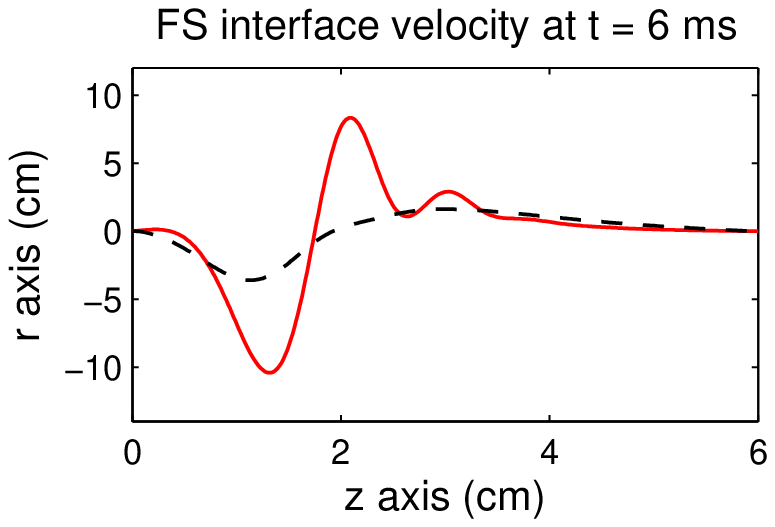}
  \includegraphics[scale=0.48]{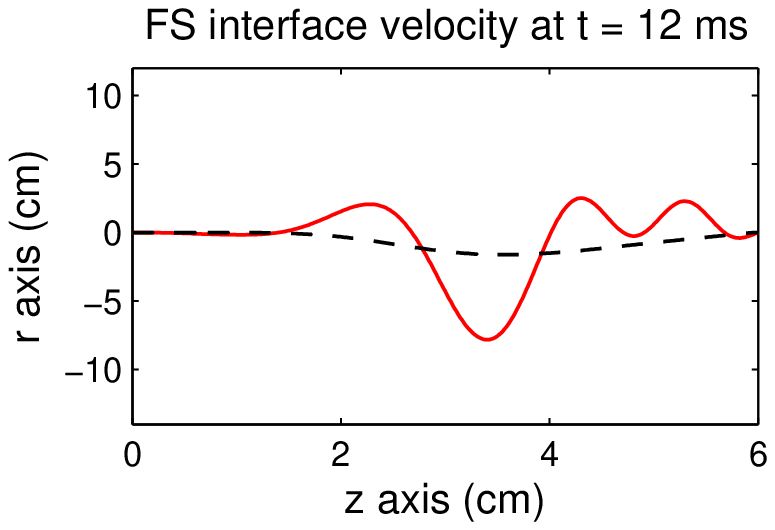}
 }
 \caption{Fluid-structure interface displacement (top) and velocity (bottom) obtained using the multilayered wall model discussed in the present manuscript, and the single layered model 
from \cite{thick}, shown at times $t = 1$ ms, $t = 6$ ms, and $t = 12$ ms.}
\label{Ex2F2}
 \end{figure}
 In the first row of Figure~\ref{Ex2F2} three snap-shots of the fluid-structure interface are shown as the inlet pressure wave travels 
 down the tube. In the second row of Figure~\ref{Ex2F2} the same three snap-shots are shown, but for the fluid-structure interface velocity.
 The red solid line in these figures corresponds to the massless fluid-structure interface in the FSI problem with a single thick structural layer.
 The black dashed line correspond to the fluid-structure interface with mass in the FSI problem with two structural layers. 
 We see significant damping of the traveling wave in the case when the fluid-structure interface has mass. This indicates that 
 inertia of the fluid-structure interface with mass regularizes solutions of FSI problems.

 This is reminiscent of the results by Hansen and Zuazua \cite{HansenZuazua} in which the presence of point mass
at the interface between two linearly elastic strings with solutions in asymmetric spaces (different
regularity on each side) allowed the proof of well-posedness due to the regularizing effects by the
point mass. 
In particular,  in \cite{HansenZuazua} two linearly elastic strings were considered, meeting at a point mass.
The elastodynamics of each string was modeled by the linear wave equation.
It was shown that as the wave with the displacement in  $H^1(0,L)$
and velocity in $L^2(0,L)$ passes through the point mass, a reflected and a transmittedf wave forms. 
The transmitted wave, which passes through the point mass, gets smoothed out to $H^2(0,L)$ regularity in displacement,
and $H^1(0,L)$ regularity in velocity. 
A numerical simulation of this phenomenon was shown in \cite{LescarretZuazua}. 
\begin{figure}[ht!]
 \centering{
 \includegraphics[scale=0.33]{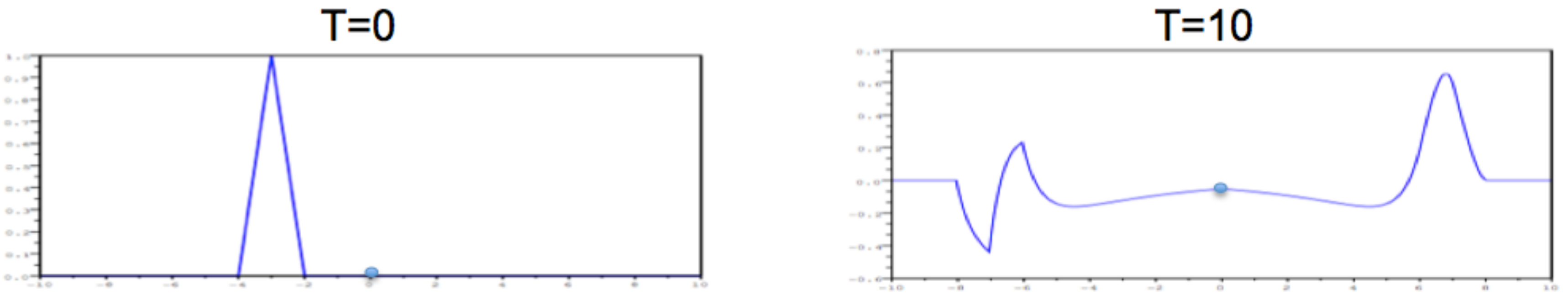}
 }
 \caption{Regularizing effects of point mass. The figure is taken from \cite{LescarretZuazua}. The initial data (left panel) is smoothed out
 as the transmitted wave traveling to the right, passes through the point mass (right panel).}
\label{fig:Zuazua}
 \end{figure}
Figure~\ref{fig:Zuazua} shows one of the results from \cite{LescarretZuazua}. The panel on the left shows the initial displacement in $H^1(0,L)$
with zero initial velocity,
located just left from the point mass. The panel on the right shows the solution at time $T = 10s$ at which the reflected and transmitted
waves have formed, with the displacement of the reflected wave on the left of the point mass still in $H^1(0,L)$,
but with the displacement of the transmitted wave, shown to the right of the point mass, belonging to $H^2(0,L)$.
 For a reader with further interest in the area of simplified coupled problems we
mention \cite{KochZauzua,RauchZhangZuazua,ZhangZuazuaARMA07}.


\if 1 = 0
This chapter addresses an operator splitting approach to study multi-physics problems related to fluid-structure interaction
with application to hemodynamics. 
The methodology is based on the Lie splitting scheme, also known as the Marchuk-Yanenko scheme. 
The splitting discussed in this chapter deals successfully with the added mass effect which is known to be responsible for
instabilities in loosely-coupled Dirichlet-Neumann schemes for FSI problems in which the density of the structure is comparable
to that of the fluid, which is the case in blood flow applications.
 Particular attention was payed to a multi-physics FSI problem in which the structure is composed of multiple 
structural layers. Problems if this kind arise, for example, in modeling blood flow through human arteries which are composed 
of several layers, each with different mechanical characteristics and thickness. A benchmark problem was studied in which 
the structure consists of two layers: 
a thin layer which is in contact with the fluid, and a thick layer which sits on top of the thin layer. The thin layer serves as
a fluid-structure interface with mass. Both analytical (existence of a weak solution) as well as numerical results were studied for the
underlying benchmark problem. In particular, it was shown that the proposed scheme converges to a weak solution to the full
nonlinear fluid-multi-layered structure interaction problem. Two academic examples were considered to test the performance of 
the numerical scheme. 

The analytical and numerical methods presented here apply with slight modifications to a larger class of problems. They include, for example,
a study of FSI with one structural layer (thin \cite{GioSun,MarSun,BukacCanicAppl}, or thick \cite{ThickPaper}), FSI with poroelastic structures \cite{Martina_Biot}, FSI between
a mechanical device called stent, arterial wall and fluid \cite{BorSunStent}, and FSI involving a non-Newtonian fluid \cite{NecasovaLukacova,LukacovaCMAME,Lukacova}.

This chapter provides the basic mathematical tools for further development of analytical and computational methods based on 
the Lie operator splitting approach, to study various multi-physics problems involving fluid-structure interaction.
\fi

\section{Appendix: Proof of convergence  as $h \to 0$ of FSI solutions with composite structures to a FSI solution with a single thick structure.}
We consider here the original FSI problem \eqref{NS1}-\eqref{dynamic} satisfying the following simplifying assumptions:
\begin{itemize}
 \item[1.] The fluid problem is defined on a fixed fluid domain $\Omega^f$ (linear coupling).
 \item[2.] Fluid advection is neglected leading to a time-dependent Stokes problem for the fluid.
\end{itemize}
In this case the coupled FSI problem with composite structure 
consisting of a thick structure sitting on top of  a thin membrane with thickness $h$, serving as a fluid-structure interface with mass,
takes the following form:
\begin{equation}\label{FSILin}
\left.
\begin{array}{rcl}
\varrho_f\partial_t{\boldsymbol v}^h &=& \nabla\cdot\boldsymbol\sigma^f({\boldsymbol v}^h,p^h) 
\\
\nabla\cdot{\boldsymbol v }^h &=&0 
\end{array}
\right\} \ {\rm in}\ {\Omega}^f\times (0,T),
\end{equation}

\begin{equation}
\rho_{s} \partial^2_t {\boldsymbol U}^h + \gamma \boldsymbol U^h = \nabla \cdot \boldsymbol \sigma^s( {\boldsymbol U}^h)\ \  \textrm{in} \; {\Omega}^s \times (0,T),
\end{equation}

\begin{equation}
\left.
\begin{array}{rcl}
\rho_{m} h \partial_t^2 \eta^h_z-C_2(h) \partial_z\eta^h_r-C_1(h) \partial^2_z\eta^h_z=-\boldsymbol\sigma^f{\bf n}\cdot{\bf e}_z+S{\bf n}\cdot{\bf e}_z
\\
\rho_{m} h \partial_t^2 \eta^h_r+C_0(h) \eta^h_r+C_2(h) \partial_z\eta^h_z=-\boldsymbol\sigma^f{\bf n}\cdot{\bf e}_r+S{\bf n}\cdot{\bf e}_r
\\
{\boldsymbol v}^h=\partial_t \boldsymbol\eta^h
\end{array}
\right\} \ 
{\rm on}\; \Gamma \times (0,T),
\end{equation}

\begin{equation}
\boldsymbol {\sigma n}^f_{in} = -p_{in}(t) \boldsymbol{n}^f_{in} \  \textrm{on} \; \Gamma^f_{in} \times (0,T), 
 \end{equation}
 \begin{equation}\label{FSILinEnd}
\boldsymbol {\sigma n}^f_{out} = -p_{out}(t) \boldsymbol{n}^f_{out} \  \textrm{on} \; \Gamma^f_{out} \times (0,T),
\end{equation}
where the coefficients $C_{i}(h), i = 0, 1, 2,$ all depend on $h$ and are given by formulas \eqref{coefficientsC}.
The super-script $h$ in ${\boldsymbol v}^h$ refers to a solution of the FSI problem with composite structure in which
$h$ denotes the thickness of the fluid-structure interface with mass.
Equations \eqref{FSILin}-\eqref{FSILinEnd}  are supplemented with initial and boundary data.
This gives rise to a well-defined {\sl linear} FSI problem. 

The weak formulation of problem~\eqref{FSILin}-\eqref{FSILinEnd} is given by:
\begin{equation}\label{WeakLin}
\begin{array}{c}
\displaystyle{
-\varrho_f\int_0^T\int_{\Omega^f}{\boldsymbol{v}}^h\cdot\partial_t\boldsymbol\varphi
+2\mu_f\int_0^T\int_{\Omega^f}{\bf D}({\boldsymbol{v}}^h):{\bf D}(\boldsymbol\varphi)
-\varrho_m h\int_0^T\int_{\Gamma}\partial_t\eta_r^h\partial_t\zeta_r
}
\\ \\
\displaystyle{
-\varrho_m h\int_0^T\int_{\Gamma}\partial_t\eta_z^h\partial_t\zeta_z
+C_2(h)\int_0^T\int_{\Gamma}\big (\partial_z\eta_z^h\zeta_r-\partial_z\eta_r^h\zeta_z \big )
+C_1(h)\int_0^T\int_{\Gamma}\partial_z\eta_z^h\partial_z\zeta_z
}
\\ \\
\displaystyle{
+C_0(h)\int_0^T\int_{\Gamma}\eta_r^h\zeta_r
-\varrho_s\int_0^T\int_{\Omega_s}\partial_t{\bf U}^h\cdot\partial_t\boldsymbol\psi
+2 \mu_{s} \int_0^T\int_{\Omega^s} \boldsymbol D(\boldsymbol U^h) : \boldsymbol D(\boldsymbol \psi) 
}
\\ \\
\displaystyle{
+\lambda_{s} \int_{\Omega^s} (\nabla \cdot \boldsymbol U^h) (\nabla \cdot \boldsymbol \psi)
+\gamma \int_0^T\int_{\Omega^s} \boldsymbol U^h \cdot \boldsymbol \psi
=\pm\int_0^T\int_{\Gamma^f_{in/out}}p_{in/out}(t)\varphi_z
}
\\ \\
\displaystyle{
+\rho_f \int_{\Omega^f}{\boldsymbol{v}}_0\cdot\boldsymbol\varphi(0)+\rho_m h \int_{\Gamma}\partial_t \boldsymbol{\eta}(0)\boldsymbol\zeta(0)
+\rho_s \int_{\Omega^s}{\bf V}_0\cdot\boldsymbol\psi(0),
}
\end{array}
\end{equation}
for all $(\boldsymbol\varphi,\boldsymbol\zeta,\boldsymbol\psi)\in C_c^1(0,T;V^f(\Omega^f)\times H_0^1(0,L)^2\times U^s(\Omega^s))$,
such that $\boldsymbol\varphi|_{\Gamma} = \boldsymbol{\zeta}=\boldsymbol\psi|_{\Gamma}$,
where
\begin{eqnarray*}
V^f &=& \{\boldsymbol\varphi \in H^1(\Omega^f)^2\ | \ \varphi_r|_{r = 0} = 0, \boldsymbol\varphi|_{z = 0,L} = 0\}, \\
U^s &=& \{\boldsymbol\varphi \in H^1(\Omega^s)^2\ | \ \psi_z|_{\Gamma_{ext}}= 0, \boldsymbol\psi|_{z = 0,L} = 0\}.
\end{eqnarray*}
and
$C_0(h), C_1(h),$ and $C_2(h)$ are given by~\eqref{coefficientsC}.

It is easy to see that for every $h>0$ there exists a weak solution $({\boldsymbol{v}}^h,\boldsymbol\eta^h,{\bf U}^h)$ 
to problem~(\ref{FSILin})-\eqref{FSILinEnd}
satisfying the following energy inequality:
\begin{equation}\label{LinEE}
\begin{array}{c}
\displaystyle{
\varrho_f\|{\boldsymbol{v}}^h\|^2_{L^{\infty}([0,T);L^2(\Omega^f))}+\mu_f\|{\bf D}({\boldsymbol{v}}^h)\|^2_{L^{2}(0,T;L^2(\Omega^f))}
+\varrho_m h\|\partial_t\boldsymbol\eta^h\|^2_{L^{\infty}(0,T;L^2(\Gamma))}
}
\\ \\
\displaystyle{h\Big (
4 \mu_{m} \|\frac{\eta^h_r}{R}\|^2_{L^2(0,L)} +
  4 \mu_{m} \|\partial_z\eta^h_z\|^2_{L^2(0,L)}
  +\frac{4 \mu_{m} \lambda_{m}}{\lambda_{m}+2\mu_{m}} 
  \|\partial_z\eta^h_z+ \frac{\eta^h_r}{R} \|^2_{L^2(0,L)}\Big )
}
\\ \\
\displaystyle{
+\varrho_s\|\partial_t{\bf U}^h\|^2_{L^{\infty}(0,T;L^2(\Omega^s))}
+C(\mu_s,\lambda_s,\gamma)\|{\bf U}^h\|^2_{L^2(0,T;H^1(\Omega^s)}
}
\\ \\
\displaystyle{
\leq C(P_{in},P_{out},{\boldsymbol{v}}_0,{\bf V}_0,\eta_0).
}
\end{array}
\end{equation}
Since the right hand side of this energy inequality does not depend on $h$, the sequence of solutions
of problems \eqref{FSILin}-\eqref{FSILinEnd}, defined for different values of $h$, is uniformly bounded.
From here, the following weak and weak* convergence results hold as $h\to 0$:
\begin{equation}\label{weak_conv}
\begin{array}{c}
{\boldsymbol{v}}^h\rightharpoonup{\boldsymbol{v}}\quad{\rm weak*}\quad{\rm in}\; L^{\infty}(0,T;L^2(\Omega^f)),
\\
{\boldsymbol{v}}^h\rightharpoonup{\boldsymbol{v}}\quad{\rm weak}\quad{\rm in}\; L^{2}(0,T;H^1(\Omega^f)),
\\
{\bf U}^h\rightharpoonup{\bf U}\quad{\rm weak}\quad{\rm in}\; L^{2}(0,T;H^1(\Omega^s)),
\\
\partial_t {\bf U}^h\rightharpoonup\partial_t{\bf U}\quad{\rm weak*}\quad{\rm in}\; L^{\infty}(0,T;L^2(\Omega^s)),
\\
\sqrt{h}\boldsymbol\eta^h\rightharpoonup \boldsymbol\eta\quad{\rm weak}\quad{\rm in}\; L^{2}(0,T;L^2(\Gamma))^2,
\\
\sqrt{ h}\ \partial_t\boldsymbol\eta^h\rightharpoonup\boldsymbol\eta_1\quad{\rm weak*}\quad{\rm in}\; L^{\infty}(0,T;L^2(\Gamma))^2,
\\
\sqrt{ h}\ \partial_z\eta_z^h\rightharpoonup \eta_2\quad{\rm weak*}\quad{\rm in}\; L^{\infty}(0,T;L^2(\Gamma)).

\end{array}
\end{equation}
From the uniqueness of the limit in ${\mathcal D}'(\Gamma)$ we conclude that $\boldsymbol\eta_1 = \partial_t \boldsymbol\eta$,
and $\eta_2 = \partial_z \eta_z$. 

What is left is to identify the limits $\boldsymbol\eta$ and $\partial_t \boldsymbol\eta$. 
From the kinematic 
coupling condition and the Trace Theorem we have:
$$
\|\partial_t\boldsymbol\eta^h\|_{L^{2}([0,T);L^2(\Gamma))}=\|{\boldsymbol{v}}^h\|_{L^{2}(0,T;L^2(\Gamma))}
\leq\|{\boldsymbol{v}}^h\|_{L^{2}(0,T;H^1(\Omega^f))}\leq C,
$$
where $C$ is independent of $h$.  From the uniform boundedness of $\|\partial_t\boldsymbol\eta^h\|_{L^{2}(0,T;L^2(\Gamma))}$
we  have
$$
\sqrt{ h}\ \partial_t\boldsymbol\eta^h\to 0\quad{\rm in}\quad L^2(0,T;L^2(\Gamma)).
$$
Uniqueness of the limit in ${\mathcal D}'(\Gamma)$ implies $\partial_t\boldsymbol\eta=0$.
Similarly, we get $\boldsymbol\eta=0$ using the fact that $\boldsymbol\eta^h(s,.)=\eta_0+\int_0^s\boldsymbol\partial_t\eta^h$. 

We can now pass to the limit as $h \to 0$ in the linear problem~\eqref{WeakLin}-\eqref{FSILinEnd} 
by using~\eqref{weak_conv} to get the following equality:
\begin{equation}\label{WeakLinThick}
\begin{array}{c}
\displaystyle{
-\varrho_f\int_0^T\int_{\Omega^f}{\boldsymbol{v}}\cdot\partial_t\boldsymbol\varphi
+2\mu_f\int_0^T\int_{\Omega^f}{\bf D}({\boldsymbol{v}}):{\bf D}(\boldsymbol\varphi)
-\varrho_s\int_0^T\int_{\Omega_s}\partial_t{\bf U}\cdot\partial_t\boldsymbol\psi
}
\\ \\
\displaystyle{
+2 \mu_{s} \int_0^T\int_{\Omega^s} \boldsymbol D(\boldsymbol U) : \boldsymbol D(\boldsymbol \psi)
+\lambda_{s} \int_{\Omega^s} (\nabla \cdot \boldsymbol U) (\nabla \cdot \boldsymbol \psi)
+\gamma \int_0^T\int_{\Omega^s} \boldsymbol U \cdot \boldsymbol \psi 
}
\\ \\
\displaystyle{
=\pm\int_0^T\int_{\Gamma^f_{in/out}}p_{in/out}(t)\varphi_z
+ \rho_f \int_{\Gamma}\boldsymbol{v}_0\boldsymbol\zeta(0)
+\rho_s \int_{\Omega^s}{\bf V}_0\cdot\boldsymbol\psi(0),
}
\end{array}
\end{equation}
for all $(\boldsymbol\phi,\boldsymbol\psi)\in C_c^1([0,T);V^f(\Omega^f)\times U^s(\Omega^s))$
such that $\boldsymbol\phi|_{\Gamma} = \boldsymbol\psi|_{\Gamma}$.
This shows that the limiting functions satisfy the weak form of the FSI with a single, thick elastic structure,
coupled to the motion of an incompressible, viscous fluid via a mass-less fluid-structure interface. 

{\bf Remark 2.}
We did not have to identify the limit $\boldsymbol\eta$ in order to pass to the limit. Namely,
we have:
$$
\varrho_m h\partial_t\boldsymbol\eta^h=\sqrt{\varrho_m h}(\sqrt{\varrho_m h}\partial_t\boldsymbol\eta^h)
\rightharpoonup\sqrt{\varrho_m h}\partial_t\boldsymbol\eta=0.
$$

\section{Acknowledgements} 
Research of Buka\v{c} was partially supported by NSF via grant DMS-1318763. Research of \v{C}ani\'{c} was partially supported
by NSF via grants NIGMS DMS-1263572, DMS-1318763, DMS-1311709, DMS-1262385 and DMS-1109189. Research of Muha
was partially supported by NSF via grant DMS-1311709. 
\bibliographystyle{plain}
\bibliography{multilayered}
\end{document}